
\magnification=1200
\hsize=14.2cm
\vsize=19.7cm
\hoffset=-0.42cm
\voffset=0.6cm

\nopagenumbers

\input amssym.def
\input amssym.tex


\def\Re{{\rm Re}\,}

\def\txt#1{{\textstyle{#1}}}
\def\scr#1{{\scriptstyle{#1}}}
\def\r#1{{\rm #1}}
\def\B#1{{\Bbb #1}}
\def\c#1{{\cal #1}}

\def\mymod#1{{\,(\bmod\,#1)}}
\def\smymod#1{{\,(\!\bmod\,#1)}}

\font\title=cmbx12 at 14pt
\font\author=cmr12
\font\srm=cmr8 at 9pt
 at11pt

\def\rightheadline{\hfil{\srm
Sieve Method and its History}
\hfil\folio}
\def\leftheadline{\rm\folio\hfil{\srm
Y. Motohashi}\hfil}
\def\emptyheadline{\hfil}
\headline{\ifnum\pageno=1 \emptyheadline\else
\ifodd\pageno \rightheadline \else \leftheadline\fi\fi}

\def\firstpage{\hss{\vbox to 1cm{\vfil\hbox{\rm\folio}}}\hss}
\def\emptyfootline{\hfil}
\footline{\ifnum\pageno=1\firstpage\else
\emptyfootline\fi}


\centerline{\title An Overview of the Sieve Method
and its History}
\vskip 0.7cm
\centerline{\author Yoichi Motohashi}
\vskip 1cm
\noindent
{\bf Preface}
\medskip
\noindent
Trying to decompose an integer into
a product of integers, we feel irritation. There should dwell the
reason why any prime appears like a real gem that one can touch and hold.
We thus muse ever and again how and when ancient people discovered 
the way of sifting out primes and began appreciating them. 
Perhaps those who conceived
the divisibility had already some sieves in their minds. Indeed,
a wealth of evidences have been excavated 
supporting our view. The story to be told below
must have originated more than five millennia ago$\,{}^{1)}$, while the primordial
intellectual irritation has remained fresh and fundamental till today.
\par
The history of the Sieve Method is rich and fascinating; we 
would need a volume to exhaust the story. In the present
article we shall instead concentrate on several pivotal ideas that
made progress possible; so the scope is inevitably limited.
Nevertheless, you will encounter instances of precious
mathematical achievements that people in the future will certainly
continue to relate.
\smallskip
Notes are to be read as essential parts, although
they are in the style of personal memoranda. 
Mathematical symbols and definitions are introduced where they are needed
for the first time, and will continue to be effective until otherwise
stated. Theorems are given somewhat implicitly, and details such as
domains of variables are to be induced from the context.
References are restricted mostly to seminal works in respective
developments. Basic facts from Analytic Number Theory could 
be found in the monographs [26] [64].
\medskip
\noindent
Remark: This is a translation of our Japanese expository article that was
published under the title `{\it An overview of sieve methods}' 
in the second issue of the 52nd volume of Sugaku, 
the Mathematical Society of Japan, April 2005. 
At this opportunity we have made some revision and changed the title
into something more appropriate. Also, it should be noted that
events in the last two decades are left untouched except for a few.
\vskip 0.7cm
\noindent
{\bf Chapter 1. Brun's Sieve}
\medskip
\noindent
{\bf 1.1} Mark any natural number that is divisible 
by the first prime 2, and repeat the same  with all 
other primes less than a given $z>2$. Then any natural number
less than $z^2$ that remains unmarked is either $1$ or a prime
in the interval $[z,z^2)$, since such an integer does not have
two or more prime factors. This is a version of the well-known sieve
method named after Eratosthenes of Alexandria$\,{}^{2)}$. 
\par
Eratosthenes' Sieve might appear to be quite effective, 
especially when $z$ is large, for it allows us to expand 
the table of all primes less than $z$ to that of those
less than $z^2$. However, if we look into the quantitative
aspect of the method or if it is required to
count the number of integers unsifted, then 
Eratosthenes' Sieve becomes virtually ineffectual$\,{}^{3)}$. 
Viggo Brun [12]
confronted the challenge of improving Eratosthenes' Sieve
to turn it into a quantitatively effective device, 
and became the founder of the modern theory of the Sieve Method$\,{}^{4)}$. 
Let us see how he brought the first light [10] into the darkness 
of 2100 years.
\medskip
\noindent
{\bf 1.2} Let $P(z)$ be the product of all primes less than $z$, and
let $(m,n)$ be the greatest common divisor of integers $m,n$. That
an integer $n$ does not have any prime factor less than $z$ is
equivalent to $(n,P(z))=1$. Thus, the characteristic function
of the set of all such integers is expressed as
$$
\sum_{d|(n,P(z))}\mu(d),\leqno(1.1)
$$
where $\mu$ is the M\"obius function. This coincides with the above
procedure of marking integers. In fact, if $n$ has $r$ marks, then
the value of the last sum is equal to $(1-1)^r$. Hence 
$(1.1)$ could be identified with Eratosthenes' Sieve$\,{}^{5)}$.
\par
We consider a finite sequence $\c{A}$ of integers, and put
$\c{S}(\c{A},z)=\{n\in{\cal A}: (n,P(z))=1\}$. Then $(1.1)$
implies that
$$
\left|\c{S}(\c{A},z)\right|=\sum_{d|P(z)}\mu(d)|\c{A}_d|,
\quad \c{A}_d=\left\{n\in\c{A}:\, n\equiv 0\mymod d\right\}.
\leqno(1.2)
$$
Thus, in order to either evaluate or
bound $\left|\c{S}(\c{A},z)\right|$,   
we need to have certain information about the behaviour of
$|\c{A}_d|$ with variable $d$; and following a general practice
we write
$$
|\c{A}_d|={\omega(d)\over d}X+R_d, \quad \omega(d)\ge0\,,\leqno(1.3)
$$
where $\omega$ is a multiplicative function. One may regard
$(\omega(d)/d)X$ as the main term and $R_d$ as the remainder;
for instance, $X$ can be seen as an approximation to $|\c{A}|$. 
The terms $R_d$ should be small either individually 
or in a certain sense of mean. Inserting $(1.3)$ into $(1.2)$,
we have
$$
\left|\c{S}(\c{A},z)\right|=V(z,\omega)X
+R({\cal A},z),\leqno(1.4)
$$
where
$$
V(z,\omega)=\prod_{p<z}\left(1-{\omega(p)\over p}\right),
\quad R({\cal A},z)
=\sum_{d|P(z)}\mu(d)R_d. \leqno(1.5)
$$
Hereafter $p$ denotes a generic prime.
\medskip
\noindent
{\bf 1.3} The identity $(1.4)$ does not have much realistic contents,
however. To see this, we consider the problem of counting primes 
in a given interval. Thus, let $x$ be sufficiently large, and put
${\cal A}=\left\{n: x-y\le n<x\right\}$, $2\le y\le x/2$. With this, we
have $\omega\equiv1$, $X=y$, $R_d=[-(x-y)/d]-[-x/d]-y/d$, where $[a]$
is the integral part of $a$. The above discussion gives
$$
\pi(x)-\pi(x-y)=V(\sqrt{x},1)y+R({\cal A},\sqrt{x}),\leqno(1.6)
$$
where $\pi(x)$ is the number of primes less than $x$ as usual. 
Since$\,{}^{6)}$
$$
V(z,1)={e^{-c_E}\over\log z}\left(1+O\left({1\over\log z}\right)\right)
\quad (\hbox{$c_E$: the Euler constant}),\leqno(1.7)
$$
we get
$$
\pi(x)-\pi(x-y)=2e^{-c_E}(1+o(1)){y\over\log x}
+R({\cal A},\sqrt{x}).\leqno(1.8)
$$
On the other hand, the Prime Number Theorem suggests that
$$
\pi(x)-\pi(x-y)=(1+o(1)){y\over\log x}.\leqno(1.9)
$$
Thus it is plausible$\,{}^{7)}$ to have $R({\cal A},\sqrt{x})
\sim(1-2e^{-c_E})(y/\log x)$. That is, each of the two terms 
on the right of $(1.6)$ can not be the main term or 
the remainder term either. It seems extremely hard$\,{}^{8)}$ to deduce 
this fact directly from the definition of $R({\cal A},\sqrt{x})$.
\medskip
\noindent
{\bf 1.4} As another example, we shall consider 
$\pi_2(x)$ the number of twin primes less than $x$. This time 
we work with the sequence ${\cal A}=\{n(n+2):1\le n<x-2\}$; in fact, 
$\left|{\cal S}({\cal A},\sqrt{x})\right|=\pi_2(x)+O(\sqrt{x})$.  
We have $\omega(2)=1$, $\omega(p)=2$ ($p\ge3$), and by $(1.4)$ and $(1.7)$
it follows, after a rearrangement, that
$$
\pi_2(x)=8e^{-2c_E}(1+o(1)){x\over(\log x)^2}\prod_{p\ge3}
\left(1-{1\over (p-1)^2}\right)+R({\cal A},\sqrt{x}).\leqno(1.10)
$$
It is, however, much more difficult to deal with this 
$R({\cal A},\sqrt{x})$ than the analogue in the last section.
In a conjecture due to Hardy and Littlewood [23] the asymptotic identity
$\pi_2(x)\sim 2Cx/(\log x)^2$ is predicted, where $C$ is the Euler product
on the right of $(1.10)$. Thus, again, each of the two terms on the right of
$(1.10)$ cannot be the main term or the remainder term either, and 
the expression $(1.10)$ does not stand for anything meaningful for the
Twin Prime Conjecture.
\medskip
\noindent
{\bf 1.5} The difficulties with $R(\c{A},z)$ observed above stems from
the fact that the number of 
summands in the defining sum $(1.5)$ may be too big to handle. As $z$
increases, the factors of $P(z)$ can become huge, and 
moreover their number too. Whether it is possible or not to
detect any dramatic cancellation among the summands should obviously be 
tremendously difficult to see. Here is the reason for
the limitation of Eratosthenes' Sieve.
\par
In 1915 Brun [10] broke this spell with a surprisingly simple idea. He
threw the explicit formula $(1.1)$ out, and replaced it by an inequality
bounding the characteristic function from above and below so that he
could gain an effective control over the size of participating
factors of $P(z)$. In other words, he moved the sieve theory from the
classic world of exactness to the modern world of reserved certainty. More
explicitly, his idea is embodied in 
$$
\sum_{\scr{d|(n,P(z))}\atop\scr{\nu(d)\le2\ell+1}}\mu(d)
\le\sum_{d|(n,P(z))}\mu(d)\le
\sum_{\scr{d|(n,P(z))}\atop\scr{\nu(d)\le2\ell}}\mu(d),
\leqno(1.11)
$$
with $\nu(d)$ the number of different prime factors of $d$; in fact
we have, for any $\ell\ge0$, 
$$
\sum_{\scr{d|m}\atop\scr{\nu(d)\le\ell}}\mu(d)=(-1)^\ell
{\nu(m)-1\choose\ell}.\leqno(1.12)
$$
The inequality $(1.11)$ is often called Brun's Pure Sieve.
\medskip
\noindent
{\bf 1.6} Let us try Brun's idea in the instance of the twin 
prime problem; we are going to bound $\pi_2(x)$ from above.
Let $z\le\sqrt{x}$ and $\ell$ be to be fixed later. We apply $(1.11)$
to the $\c{A}$ of Section 1.4, and have
$$
\eqalign{
\pi_2(x)&\le |{\cal S}({\cal A},z)|+z\cr
&\le x\sum_{\scr{d|P(z)}\atop\scr{\nu(d)\le2\ell}}
{\mu(d)\omega(d)\over d}+\sum_{d\le z^{2\ell}}\omega(d)+z,
}\leqno(1.13)
$$
in which we have used $|R_d|\le\omega(d)$. Compared with the sum of
divisors, the second sum over $d$ is $O(\ell z^{2\ell}\log z)$. The
difference between the first sum and $V(z,\omega)$ is
$$
-\sum_{\scr{d|P(z)}\atop\scr{\nu(d)\ge2\ell+1}}
{\mu(d)\omega(d)\over d}
\ll 2^{-2\ell}\sum_{d|P(z)}
{2^{2\nu(d)}\over d};\leqno(1.14)
$$
the symbol $\ll$ indicates in general that the 
absolute value of the left
side is less than a constant multiple of the right side. The last sum is
$O\left(V(z,1)^{-4}\right)$. Collecting these and setting 
$z=\exp(\log x/(100\log\log x))$, $\ell=\left[\log x/(4\log z)\right]$,
we obtain$\,{}^{9)}$
$$
\pi_2(x)\ll x\left({\log\log x\over\log x}\right)^2
\quad\hbox{\rm or}
\quad {\pi_2(x)\over\pi(x)}\ll {(\log\log x)^2\over\log x}.
\leqno(1.15)
$$
Therefore, twin primes occur far less frequently than ordinary primes. 
It is literally hopeless to deduce this fact via $(1.10)$. It is 
amazing that the imperfect $(1.11)$ could yield any result that
the exact $(1.1)$ would never be capable of$\,{}^{10)}$.
\medskip
\noindent
{\bf 1.7} Brun's Sieve, the title of the present chapter,
is an improved version [12] of his Pure Sieve; a combinatorial
sophistication was introduced into the choice of divisors in
$(1.11)$.  It enabled Brun to achieve the impressive bound
$$
\pi_2(x)\ll {x\over (\log x)^2}.\leqno(1.16)
$$ 
In view of Hardy--Littlewood's conjecture mentioned above, this 
should be the best possible, save for the implied constant. 
The construction of Brun's Sieve is, however, so intricate that we have to
skip the details$\,{}^{11)}$, and state only the conclusion $(1.22)$ below, 
moreover in a considerably abridged fashion. Nevertheless, we shall
reach an equivalent result in Section 3.5 via a somewhat 
different argument.
\par
For the present and later purpose, we need to rearrange the specifications
introduced in [12], to be in accordance with today's practice. 
Thus, it is customary to suppose that it holds that for any $2<z_1<z_2$
$$
{V(z_1,\omega)\over V(z_2,\omega)}=\prod_{z_1\le p<z_2}
\left(1-{\omega(p)\over p}\right)^{-1}
=\left({\log z_2\over\log z_1}\right)^\kappa
\left(1+O\left({1\over\log z_1}\right)\right),
\leqno(1.17)
$$
with a $\kappa>0$, the dimension of the sieve problem
under consideration$\,{}^{12)}$. This is to be compared with
$(1.7)$. For example, in Section 1.4 the Twin Prime Conjecture
is discussed as a sieve problem of dimension 2. On the other hand,
the inequality $(1.11)$ is understood as a particular
construction of sieve weights $\rho_r(d)$ such 
that for any integer $n$
$$
(-1)^r\left\{\sum_{d|(n,P(z))}\mu(d)-
\sum_{d|(n,P(z))}\mu(d)\rho_r(d)\right\}\ge0,
\quad \rho_r(1)=1,
\leqno(1.18)
$$
with the convention $\rho_{r_1}\equiv\rho_{r_2}$ 
($r_1\equiv r_2\mymod 2$). We have, because of $(1.3)$, 
$$
(-1)^r\left\{|{\cal S}({\cal A},z)|- V(z,\omega\,;
\rho_r)X\right\}\ge(-1)^r R({\cal A},z\,;\rho_r),\leqno(1.19)
$$
where
$$
V(z,\omega\,;\rho_r)=\sum_{d|P(z)}\mu(d)\rho_r(d)
{\omega(d)\over d},\quad
R({\cal A},z\,;\rho_r)=\sum_{d|P(z)}\mu(d)\rho_r(d)R_d.
\leqno(1.20)
$$
According as $r\equiv 0,\,1\mymod 2$, the inequality $(1.19)$ is called
the lower and the upper bounds of $|{\cal S}({\cal
A},z)|$; occasionally only one of them is considered.
Naturally, we shall regard $V(z,\omega\,;\rho_r)X$ as the main term, and
$R({\cal A},z\,;\rho_r)$ as the remainder. In order to have any
effective control over the latter, we ought to  impose a limitation to the
size of participating $d$. For this sake we introduce the condition
$$
\rho_r(d)=0,\quad d\ge D,\quad d|P(z).\leqno(1.21)
$$
The parameter $D$ is called the level of the sieve weights
$\rho_r$. A sieve problem is to find $\rho_r$ that yield 
any good main term under these specifications$\,{}^{13)}$.
\par
We now exhibit the principal result of Brun [12], with a drastic
simplification: Let $D=z^\tau$. Then there exist characteristic 
functions $\rho_r$ such that
$$
V(z,\omega\,;\rho_r)=(1+O(e^{-{1\over2}\tau\log\tau}))
V(z,\omega).\leqno(1.22)
$$
\medskip
\noindent
{\bf 1.8} Let us apply the last assertion to the situation
treated in Section 1.4. We set $z^\tau=\sqrt{x}$; then the remainder term
can obviously be ignored. With a sufficiently large $\tau$, we
get the assertion
$$
\left|\left\{n<x: \hbox{if $p|n(n+2)$ then $p\le x^{1/(2\tau)}$}
\right\}\right|
\asymp {x\over(\log x)^2}.\leqno(1.23)
$$
This implies not only $(1.16)$ but also that there exist infinitely
many pairs $(n,n+2)$ of integers such that the number of prime factors
of each is less than $2\tau$.
\par
In this way, Brun [12] accomplished the first definitive advance toward the
Twin Prime Conjecture. 
If we use instead his Pure Sieve, then it
would be required to set $\ell=[25\log\log x]$ and thus $\tau\approx
\log\log x$.  Hence Brun's Sieve is far stronger. Certainly the same
could be asserted about Goldbach's Conjecture$\,{}^{14)}$; we need only to
move to the sequence $\{n(N-n): 3\le n\le N-3\}$ with an even 
integer $N$.
\medskip
\noindent
{\bf 1.9} It is possible to treat the Twin Prime Conjecture as 
a one dimensional or a linear sieve problem. To see this, let $\varphi$ be
the Euler totient function and ${\rm li}$ the logarithmic integral;
and put
$$
\pi(x;k,l)=|\{p<x: p\equiv l\mymod k\}|,
\quad E(x;k,l)=\pi(x;k,l)-{1\over\varphi(k)}{\rm li}\, x.\leqno(1.24)
$$ 
The latter is called the remainder term in the prime number theorem
for arithmetic progressions. To the sequence 
${\cal A}=\{p+2:\,3\le p<x\}$, we apply Brun's Sieve; thus
$X={\rm li}\,x$, $\omega(2)=0$, $\omega(p)=p/(p-1)$, $(p\ge3)$, and
$\kappa=1$ as well as $R_d=E(x;d,-2)$. Let $\tau$ be sufficiently large
as before. Then we have
$$
|{\cal S}({\cal A},z)|=(1+O(e^{-{1\over2}\tau\log\tau}))
V(z,\omega)\cdot {\rm li}\, x
+O\left(\sum_{\scr{d<z^\tau}\atop
\scr{2\nmid d}}|E(x;d,-2)|\right).\leqno(1.25)
$$
This time, the issue is how to deal with the second term
on the right side. We assume that
$$
\sum_{q<Q}\,\max_{\scr{a\, (\!\bmod q)}\atop
\scr{(a,q)=1}}|E(x;q,a)|\ll {x\over(\log x)^A},
\quad A\ge 3,\leqno(1.26)
$$
Then we have
$$
\left|\left\{p< x: \hbox{$p+2$ does not have prime
divisors less than $Q^{1/\tau}$}\right\}\right|\asymp{x\over(\log Q)(\log
x)}.\leqno(1.27)
$$
For example, the Extended Riemann Hypothesis would allow us to set
$Q=\sqrt{x}/(\log x)^{A+1}$; and  
$$
\left|\left\{p:\, \nu(p+2)\le 2\tau+1\right\}\right|=\infty.\leqno(1.28)
$$
Goldbach's Conjecture could be treated in much the same way. 
Comparing this with the assertions in the previous section, we apprehend 
the strength of the Extended Riemann Hypothesis$\,{}^{15)}$.
\smallskip
At $(1.25)$--$(1.26)$ is observed a typical instance of 
the relation between a sieve problem and the distribution 
of the relevant sequence of integers among arithmetic progressions.
That is, a sieve problem upon a particular sequence of integers
is reduced to the discussion of the distribution of the elements of the
sequence among arithmetic progressions with variable moduli. One may see
that  if the decomposition into arithmetic progressions is made too fine,
then relevant moduli could become too large to handle the remainder term in
the sieve  effectively. Also it might be conceivable that relatively small
moduli could contribute substantially. This situation reminds us of 
the Circle Method of Ramanujan and Hardy$\,{}^{16)}$; 
that is, there appears to
exist a relation between sieve problems and the Farey sequence. We shall
encounter the same in the next section but with a somewhat different
context.
\medskip
\noindent
{\bf 1.10} Brun's Sieve yields remarkable assertions not only
about those great conjectures but also about fundamental
queries in the theory of the distribution of primes:
$$
\leqalignno{
&\pi(x)-\pi(x-y)\ll{y\over\log y}\quad (2\le y\le x-2),&(1.29)\cr
&\pi(x;k,l)\ll{x\over \varphi(k)\log(x/k)}\quad (2\le k\le x/2).
&(1.30)
}
$$
In fact, $(1.29)$ follows immediately from an obvious modification of 
the discussion in Section 1.3. 
On the other hand, as to $(1.30)$ we may certainly assume that $(k,l)=1$; 
then for the sequence ${\cal A}=\{n< x: n\equiv l\mymod k\}$
the specification $(1.3)$ holds with $X=y/k$, $|R_d|\le1$; $\omega(p)=0$
as $p|k$, and $\omega(p)=1$ as $(p,k)=1$. In particular, we have
$V(z,\omega)\le (k/\varphi(k))V(z,1)$. The rest is straightforward. The
assertion $(1.30)$ is called traditionally the Brun--Titchmarsh theorem.
\par
As a matter of fact, if the uniformity as displayed 
prominently in $(1.29)$--$(1.30)$
is required, any currently available analytic method 
which relies on the theory of the Riemann zeta and the Dirichlet
$L$-functions is unable to produce anything comparable to the last two
assertions. Even under the Extended Riemann Hypothesis, the situation would not
change. It is indeed amazing that such an elementary idea as Brun's could 
ever take us close to the very subtlety of the distribution of primes. 
\vskip 0.7cm
\noindent
{\bf Chapter 2. Linnik's and Selberg's Sieves}
\medskip
\noindent
{\bf 2.1} Some 20 years had passed since Brun's fundamental work [12]
when Ju.\ V. Linnik [41] marked a new departure in the Sieve Method; and in a
few years A. Selberg [68]--[70] made an independent leap. Selberg's idea
was a distinctive incision into the sieve theory general. The construction
of his sieve, i.e., his sieve weights, is fundamentally different from
Brun's; thus he brought a structural change into the Sieve Method. 
On the other hand, Linnik's idea would later be
appreciated because not only of its highly 
effective sieve effect but also of the argument itself that
he employed. With Linnik's seminal work, 
a general principle was born, which
has been and is still a driving force behind many of major works in
Analytic Number Theory. This is now called the Large Sieve, 
in a much wider context than the title of his work indicated then. 
\par
We shall show the essentials of the two ideas. Also, a duality relation
between them will be disclosed. Further, it will be witnessed that
Large Sieve yields not only a sieve bound but also a spectacular
conclusion about the distribution of primes in arithmetic progressions.
\medskip
\noindent
{\bf 2.2} We need first to change the technical definition 
introduced in Section 1.7
of a sieve problem into a conceptual setting. Thus, let  $\Omega(p)$ 
be a set of residue classes $\bmod\,p$, and let us
write $n\in\Omega(p)$ in stead of $n\mymod p\in\Omega(p)$. We put
$$
{\cal S}({\cal A}, z\,;\Omega)=\left\{n\in{\cal A}: n \not\in
\Omega(p), \forall p<z\right\},\leqno(2.1)
$$
with any sequence $\c{A}$ of integers$\,{}^{17)}$. A sieve problem is newly
defined to be the estimation of 
$\left|{\cal S}({\cal A},z\,; \Omega)\right|$,
though we shall consider mainly the situation where
$\c{A}$ is the interval
$$
{\cal N}=[M,M+N)\cap{\Bbb Z},
\quad M\in{\Bbb Z},\,N\in{\Bbb N}\,.\leqno(2.2)
$$
\par
In Section 1.4, the Twin Prime Conjecture is treated with
$M=0$, $N=[x-2]$, $\Omega(2)=\{0\}$, $\Omega(p)=\{0,-2\}$
$(p\ge3)$, $|\Omega(p)|=\omega(p)$. In this example, we have 
$|\Omega(p)|\le2$; but in general we should not restrict
the size of $\Omega(p)$. However, if we have, for instance, always
$|\Omega(p)|\ge cp$ with a certain fixed $c>0$, Brun's Sieve
is not applicable, for the condition $(1.17)$ on the sieve dimension
is violated. Then, is there any sieve procedure that is effective
even when $|\Omega(p)|$ may become huge like that? It was Linnik
[41] who gave the first answer to this intriguing problem. 
\medskip
\noindent
{\bf 2.3} Linnik argued as follows: Let $\varpi$ be the
characteristic function of the set $\{n\in{\Bbb Z}: n
\not\in\Omega(p), \forall p<z\}$, and put
$$
\eqalign{U(\theta)&=\sum_{M\le n< M+N}
\varpi(n)\exp(2\pi i n\theta),\cr
U(\theta;p,a)&=\sum_{\scr{M\le n< M+N}\atop
\scr{n\equiv a\smymod p}}\varpi(n)\exp(2\pi in\theta).
}\leqno(2.3)
$$ 
Then we have
$$
\sum_{a=1}^{p-1}\left|U\left(\theta+{a\over p}\right)\right|^2
=p\sum_{a=1}^p|U(\theta;p,a)|^2-|U(\theta)|^2.
\leqno(2.4)
$$
On the other hand, since $U(\theta;p,a)=0$ if $a\in\Omega(p)$, $p<z$,
we have also
$$
|U(\theta)|^2=\left|\sum_{a=1}^pU(\theta;p,a)\right|^2\le (p-|\Omega(p)|)
\sum_{a=1}^p|U(\theta;p,a)|^2,\quad p<z,\leqno(2.5)
$$
or$\,{}^{18)}$
$$
|U(\theta)|^2{|\Omega(p)|\over p-|\Omega(p)|}
\le \sum_{a=1}^{p-1}\left|U\left(\theta+{a\over p}\right)\right|^2,
\quad p< z.\leqno(2.6)
$$
By the decomposition law of residue classes, we have
$$
(U(0))^2\prod_{p|q}{|\Omega(p)|\over p-|\Omega(p)|}
\le \sum_{\scr{a \smymod q}\atop\scr{(a,q)=1}}
\left|U\left({a\over q}\right)\right|^2,\quad q|P(z).\leqno(2.7)
$$
In this way, we are led to$\,{}^{19)}$
$$
\left|{\cal S}({\cal N},z\,;\Omega)\right|^2G(z,\Omega)
\le \sum_{\scr{q<z}}\sum_{\scr{a \smymod
q}\atop\scr{(a,q)=1}}\left|\sum_{M\le n<M+N}\varpi(n)
\exp\left(2\pi i {a\over q}n\right)\right|^2,\leqno(2.8)
$$
where
$$
G(z,\Omega)=\sum_{q<z}\mu(q)^2H(q,\Omega),
\quad H(q,\Omega)=\prod_{p|q}{|\Omega(p)|\over p-|\Omega(p)|}\,.
\leqno(2.9)
$$
\par
Hence, our initial problem has been reduced to the estimation of
the sum in $(2.8)$ over the Farey sequence. We skip the relevant 
discussion by Linnik himself, since we shall show an assertion
better than his, in Section 2.6. What is essential here is to
follow the procedure leading to $(2.8)$, which is termed 
Linnik's Sieve$\,{}^{20)}$. Nevertheless, the final assertion should be
displayed: By the inequality $(2.30)$, we have
$$
\left|{\cal S}({\cal N},z\, ;\Omega)\right|\le{1\over G(z,\Omega)}
(N+z^2).\leqno(2.10)
$$
\medskip
\noindent
{\bf 2.4} Selberg argued as follows: We extend $\Omega$ multiplicatively,
and let $\lambda$ be an arbitrary real-valued function such that
$\lambda(1)=1$. We have, for any integer $n$,

$$
\varpi(n)\le\left(\sum_{\scr{n\in\Omega(d)}\atop
\scr{d|P(z)}}\lambda(d)\right)^2=
\sum_{\scr{n\in\Omega([d_1,d_2])}\atop
\scr{d_1,d_2|P(z)}}
\lambda(d_1)\lambda(d_2),\leqno(2.11)
$$
where $[d_1,d_2]$ is the least common multiple of $d_1,\,d_2$. 
This inequality is trivial; nevertheless, its consequence is
impressive. 
\par
For the sake of simplicity, we assume that $\lambda(d)=0$ 
either if $d\ge z$ or if $\mu(d)=0$. Summing $(2.11)$ over $n\in\c{N}$ and
exchanging the order of summation, we have
$$
|{\cal S}({\cal N},z\,;\Omega)|\le N \cdot S+R,\leqno(2.12)
$$
where
$$
\eqalign{
&S=\sum_{d_1,d_2<z}{|\Omega([d_1,d_2])|\over
[d_1,d_2]}\lambda(d_1)\lambda(d_2),\cr 
&|R|\le \sum_{d_1,d_2<z}|\Omega([d_1,d_2])|
|\lambda(d_1)\lambda(d_2)|.
}\leqno(2.13)
$$
Selberg computed the minimum of the quadratic form $S$ with the
side condition $\lambda(1)=1$. His reasoning is illuminating$\,{}^{21)}$.
We note first
$$
{|\Omega([d_1,d_2])|\over[d_1,d_2]}
={|\Omega(d_1)|\over d_1}\cdot{|\Omega(d_2)|\over d_2}
\cdot {(d_1,d_2)\over|\Omega((d_1,d_2))|},\leqno(2.14)
$$
and by the M\"obius inversion 
$$
{(d_1,d_2)\over|\Omega((d_1,d_2))|}
=\sum_{\scr{f|d_1}\atop\scr{f|d_2}}{1\over|\Omega(f)|}\prod_{p|f}
\left(p-|\Omega(p)|\right),
\quad \mu(d_1)\mu(d_2)\ne0.\leqno(2.15)
$$
Thus
$$
S=\sum_{f<z}{\mu(f)^2\over|\Omega(f)|}\prod_{p|f}
(p-|\Omega(p)|)\cdot\xi(f)^2,\quad
\xi(f)=\sum_{\scr{d<z}\atop\scr{d\equiv0\mymod f}}
{|\Omega(d)|\over d}\lambda(d).\leqno(2.16)
$$
The inversion of the linear transform
$\lambda\mapsto\xi$ is given by
$$
\lambda(d)={d\over|\Omega(d)|}\sum_{g<z/d} \mu(g)\xi(dg).
\leqno(2.17)
$$
The side condition $\lambda(1)=1$ is thus transformed accordingly,
and
$$
S={1\over G(z,\Omega)}+\sum_{f<z}{\mu(f)^2\over|\Omega(f)|}\prod_{p|f}
(p-|\Omega(p)|)\left(\xi(f)-{1\over G(z,\Omega)}
\mu(f)H(f,\Omega)\right)^2.\leqno(2.18)
$$
Hence the optimal $\xi$ is found, and inserting it into $(2.17)$
we are led to the assertion that
$$
S={1\over G(z,\Omega)},\quad\lambda(d)={\mu(d)
\over G(z,\Omega)}\prod_{p|d}{p\over
p-|\Omega(p)|}\sum_{\scr{g<z/d}\atop\scr{(d,g)=1}}
\mu(g)^2H(g,\Omega).\leqno(2.19)
$$
Obviously the assumption on $\lambda$ imposed above is satisfied
by this specialization. Moreover, we have
$$
|\lambda(d)|\le\mu(d)^2.\leqno(2.20)
$$
In fact, we have, for any $d<z$ $(\mu(d)\ne0)$, 
$$
G(z,\Omega)=\sum_{f|d}\mu(f)^2H(f,\Omega)
\sum_{\scr{g<z/f}\atop\scr{(d,g)=1}}
\mu(g)^2H(g,\Omega),\leqno(2.21)
$$
from which $(2.20)$ follows immediately. Collecting these, we find that
$$
\left|{\cal S}({\cal N},z\,;\Omega)\right|
\le {N\over G(z,\Omega)}+R,\quad
|R|\le\left(\sum_{d<z}|\Omega(d)|\right)^2.\leqno(2.22)
$$
The procedure of the present section is called Selberg's Sieve$\,{}^{22)}$.
\medskip
\noindent
{\bf 2.5} A comparison of $(2.22)$ with $(2.10)$ might
cause the incorrect impression that Selberg's Sieve is inferior to
Linnik's$\,{}^{23)}$. In fact, $(2.10)$ 
could be deduced with Selberg's Sieve as well$\,{}^{24)}$.
\par
In order to show this, we express the characteristic function of the set
$\{n\in{\Bbb Z}: n\in\Omega(d)\}$ as
$$
\eqalign{
{1\over d}&\sum_{a\smymod d}\,\sum_{h\in\Omega(d)}\exp
\left(2\pi i(n-h){a\over d}\right)\cr
={1\over d}&\sum_{q|d}\sum_{\scr{a\smymod q}\atop\scr{(a,q)=1}}
\left(\sum_{h\in\Omega(d)}
\exp\left(-2\pi i{a\over q}h\right)\right)\cdot\exp\left(2\pi i{a\over
q}n\right).}\leqno(2.23)
$$
Inserting this into $(2.11)$, we have
$$
\left|{\cal S}({\cal N},z\,;\Omega)\right|\le
\sum_{M\le n< M+N}\left|\sum_{q<z}
\sum_{\scr{a\smymod q}\atop\scr{(a,q)=1}}b
\left({a\over q}\right)\exp\left(2\pi i{a\over q}n\right)
\right|^2,\leqno(2.24)
$$
where
$$
b\left({a\over q}\right)=\sum_{\scr{d< z}
\atop\scr{d\equiv0\smymod q}}{\lambda(d)\over d}
\sum_{h\in\Omega(d)}\exp\left(-2\pi i{a\over q}h\right).
\leqno(2.25)
$$
Applying the inequality $(2.31)$ below to the right side of $(2.24)$,
we get
$$
\left|{\cal S}({\cal N},z\,;\Omega)\right|\le (N+z^2)
\sum_{q<z}
\sum_{\scr{a\smymod q}\atop\scr{(a,q)=1}}
\left|b\left({a\over q}\right)\right|^2.\leqno(2.26)
$$
The double sum is equal to
$$
\sum_{d_1,d_2<z}{\lambda(d_1)\lambda(d_2)\over d_1d_2}
\sum_{h_1\in\Omega(d_1)}
\sum_{h_2\in\Omega(d_2)}
\sum_{q|(d_1,d_2)}\sum_{\scr{a\smymod q}\atop\scr{(a,q)=1}}
\exp\left(2\pi i{a\over q}(h_1-h_2)\right),\leqno(2.27)
$$
which coincides with $S$ above, as can readily be seen by 
observing the multiplicativity of the construction. Hence we have
$$
\left|{\cal S}({\cal N},z\,;\Omega)\right|\le (N+z^2)\cdot S.\leqno(2.28)
$$
By the argument of the previous section, up to $(2.19)$, we obtain $(2.10)$
again.
\par
Therefore the important upper bound $(2.10)$ has been proved in two ways. 
There is an obvious duality between them. It should be interesting to know
that there is such an intrinsic relation between Linnik's and Selberg's
ideas which occurred independently. By the way, the inequalities $(2.8)$
and $(2.24)$ are typical instances of applications of the Large Sieve.
\medskip
\noindent
{\bf 2.6}
We now exhibit the fundamental inequality of the Large Sieve: Let
$\{\psi_m\}$ be a finite set in a Hilbert space equipped with
the inner product $\langle\,,\rangle$. Then we have, for any $\psi$ 
in the space,${}^{25)}$
$$
\sum_m{|\langle \psi,\psi_m\rangle|^2\over
\sum_n|\langle \psi_m,\psi_n\rangle|}\le \langle 
\psi,\psi\rangle.\leqno(2.29)
$$
From this, a set of useful inequalities follow. Among them, the
following two are utilised in the above: Let $\{\theta_r\}$ be
a sequence in the unit interval, whose elements are well separated
with the minimum distance $\delta>0$ $\mymod 1$. Then we have,
for any complex vectors $\{a_n\}$, $\{b_r\}$,
$$
\sum_r \left|\sum_{M\le n<M+N}a_n\exp(2\pi in\theta_r)\right|^2
\le(N-1+\delta^{-1})\sum_{M\le n< M+N}|a_n|^2\leqno(2.30)
$$
and
$$
\sum_{M\le n< M+N}\left|\sum_r b_r\exp(2\pi in\theta_r)\right|^2
\le(N-1+\delta^{-1})\sum_r|b_r|^2,\leqno(2.31)
$$
where the interval $[M,M+N)$ is as in $(2.2)$.
The latter is a consequence of the former, and vice versa.
This is due to the well-known fact that
the norms of a bounded linear operator and its adjoint
acting in a Hilbert space are equal to each other.
\medskip
\noindent
{\bf 2.7} Hence, as far as finite intervals are concerned, Linnik's and
Selberg's Sieves give rise to the same upper bound $(2.10)$. Here are
a few assertions that are consequences of $(2.10)$:
$$
\leqalignno{
&\pi(x)-\pi(x-y)\le 2(1+o(1)){y\over\log y}\,,&(2.32)\cr
&\pi(x;k,l)\le 2(1+o(1))
{x\over\varphi(k)\log{(x/k)}}\,,
&(2.33)\cr
&\pi_2(x)\le 16(1+o(1))
\prod_{p\ge3}\left(1-{1\over (p-1)^2}\right)
{x\over (\log x)^2}\,.&(2.34)
}
$$
These can be proved, via $(2.10)$, with the corresponding
$\cal N$ and $\Omega$, with $z=(N/\log N)^{1/2}$. The asymptotic
evaluation of $G(z,\Omega)$ should not cause any difficulty.
Another interesting application could be obtained with
$\Omega(p)$ being the set of all quadratic 
non-residues $\mymod p$. 
\par
Thus, it is understood that as far as upper bounds are
concerned Linnik's and Selberg's Sieves are superior to
Brun's$\,{}^{26)}$. For instance, the bound $(2.34)$ should be
compared with the Hardy--Littlewood conjecture mentioned above.
Also, the new form $(2.33)$ of the Brun-Titchmarsh theorem
draws special attention. This is because of the following fact:
If there exist two absolute constants $\alpha,\beta>0$ with
which we have, uniformly for $(k,l)=1$,
$$
\pi(x;k,l)\le2(1-\alpha){x\over\varphi(k)\log (x/k)},
\quad k< x^\beta,\leqno(2.35)
$$
then Dirichlet $L$-functions $L(s,\chi)$, $\chi\mymod k$, should not
have any exceptional zero; so the theory of the distribution
of primes in arithmetic progressions would 
fundamentally be improved$\,{}^{27)}$. 
Also, an effective lower bound, which is
essentially best possible,  would follow for the class numbers of imaginary
quadratic number fields. In this context, it should  
be noted specifically that the critical bound
$$
\pi(x;k,l)\le 2{x\over\varphi(k)\log{(x/k)}}\leqno(2.36)
$$
has been established via a more careful application of Linnik's
Sieve$\,{}^{28)}$.
\medskip
\noindent
{\bf 2.8} The above discussion might be termed as an account of
the additive Large Sieve, for it concerns additive characters as
is indicated by $(2.30)$--$(2.31)$. We have seen the
appearance of important upper bounds in the prime number theory. In the
present section we turn to an account of the multiplicative 
Large Sieve, concerning instead Dirichlet characters; and we shall see
that  there emerges a surprising assertion on the asymptotic theory of
the distribution of primes in arithmetic progressions. More precisely,
the multiplicative Large Sieve opens a way to avoid the Extended Riemann
Hypothesis$\,{}^{29)}$. This fascinating theory was inseminated by A.
R\'enyi [65]$\,{}^{30)}$.  In order to appreciate his contribution, 
we need to review briefly the history of the theory of the distribution 
of primes.
\smallskip
\noindent
{\it Primes in Short Intervals\/}: Under the Riemann Hypothesis, the asymptotic
formula $(1.9)$ holds with, e.g., $x^{1/2}(\log x)^3<y<x/2$. 
However, G. Hoheisel [24] established, without any hypothesis, the
asymptotic formula upon the condition $x^{\vartheta}<y<x/2$, where
$\vartheta$ is a positive absolute constant less than $1$. That
was the unprecedented event in which was discovered the possibility to
avoid  the Riemann Hypothesis; and it was the beginning of the modern theory
of the distribution of primes. At the core of Hoheisel's argument
is a statistical study of the distribution of the complex zeros of the
Riemann zeta-function $\zeta(s)$ or a statistical proof of the 
Riemann hypothesis$\,{}^{31)}$, due to H. Bohr and E. Landau [3]. In
Riemann's explicit formula for the function $\pi(x)$ there is a sum over
the complex zeros, into which Hoheisel introduced the statistical study,
along with a certain zero-free region on the left of the vertical 
line $\Re s=1$. What should not be missed to observe in Bohr--Landau's
theory, especially in the context of our present discussion, 
is the r\^ole of a version of the mean values of
$\zeta(s)$. Taking later developments into account, this concerns
the analysis of
$$
\int_{-T}^T\left|\zeta\left(\txt{1\over2}+it\right)\right|^2
\Big|\sum_{n\le N}a_n n^{it}\Big|^2dt,\leqno(2.37)
$$
with $N, T\ge1$ and complex $a_n$ which are to be chosen 
appropriately$\,{}^{32)}$.
\smallskip
\noindent
{\it Least Prime Number Theorem\/}: If one wants to establish an analogue
of Hoheisel's assertion for primes in arithmetic progressions, then the
study of $\pi(x;k,l)$ $((k,l)=1)$ should be developed on the supposition
$x^\vartheta<x/k$ with a new positive absolute constant $\vartheta$ less
than $1$. This is to look for a way to avoid the Extended Riemann
Hypothesis. Obviously an extension of $(2.37)$ to Dirichlet $L$-functions
$L(s,\chi)$ $(\chi\mymod k)$ is required; but this part of the theory
does not cause any essential difficulty, for it suffices to exploit the
orthogonality of the characters. However, the theory of the distribution
of zeros of Dirichlet $L$-functions lacks what corresponds to the
zero-free region of $\zeta(s)$ mentioned above. One has to find a way to
negate this defect, which certainly requires to develop the statistical
study of the zeros in a far refined fashion than 
the followers of Bohr and Landau did. It was Linnik [42, I] who 
overcame this genuine difficulty$\,{}^{33)}$. 
There a definitive r\^ole was
played by the Brun--Titchmarsh theorem $(1.30)$$\,{}^{34)}$. Further, the
possibility of exceptional zeros caused another difficulty, or a
quantitative study of the Deuring--Heilbronn theory 
had to be developed. That was achieved by Linnik in [42, II]$\,{}^{35)}$. 
In this way the Least Prime Number Theorem was established; that is, there exists an
absolute constant $c>0$ such that the 
least prime in every reduced residue class $\bmod\, k$ is less than $k^c$.
\smallskip
\noindent
{\it Mean Prime Number Theorem\/}: Thus Linnik found a way to avoid 
the Extended Riemann Hypothesis. It concerned, however, 
a single modulus, though the uniformity on it was of course maintained.
Thus the next target was to find a way to avoid the Extended Riemann Hypothesis
simultaneously for all moduli in an arbitrary finite range. This time,
a genuine difficulty took place in extending $(2.37)$, which is the
principal difference from Linnik's situation. It was R\'enyi [65] who
resolved this difficulty. He started with Linnik's fundamental work [41],
and developed a version of the multiplicative Large Sieve to extend
$(2.37)$ to a double sum over moduli and characters, analogously involving
Dirichlet $L$-functions and polynomials twisted by characters. 
\par 
Then, he could establish $(1.26)$ for $Q=x^\alpha$ with an
absolute constant $\alpha>0$, without any hypothesis; this is
R\'enyi's Mean Prime Number Theorem$\,{}^{36)}$. Hence, as can
be seen from $(1.25)$--$(1.27)$, R\'enyi superseded Brun and made a
great step toward the Twin Prime and the Goldbach Conjectures.
\par
Naturally, efforts afterward were  
concentrated on the improvement$\,{}^{37)}$ of
R\'enyi's theorem; that is, to find larger
$Q$. Finally, after 17 years of a series of struggles, it was
established that
$$
\hbox{The inequality $(1.26)$ holds with
$\displaystyle {Q={\sqrt{x}\phantom{{}^B}
\over(\log x)^B}}$},
\leqno(2.38)
$$
where $B$ is a function of $A$. This is called  
E. Bombieri--A.I. Vinogradov's Mean Prime Number Theorem$\,{}^{38)}$.
Bombieri's argument [4] stands on the tradition of the Large Sieve;
and Vinogradov [77] relied on the Dispersion Method [43], 
another fundamental invention of Linnik$\,{}^{39)}$. Despite the difference
in their methods, what they achieved is essentially equivalent
to each other and to the consequence of the Extended Riemann Hypothesis,
especially in the context of its applications 
to sieve problems as exhibited above. In Bombieri's argument$\,{}^{40)}$,
$(2.38)$ could be said to be a consequence of the inequality
$$
\sum_{q< Q}{q\over\varphi(q)}
\,\mathop{\sum\!\raise 3pt\hbox {${}^*$}}_{\chi\smymod q}\left|
\sum_{M\le n< M+N}a_n\chi(n)\right|^2\le
(N-1+Q^2)\sum_{M\le n< M+N}
|a_n|^2,\leqno(2.39)
$$
where the asterisk means that the sum is restricted to primitive
characters. Connecting multiplicative characters with additive
characters via Gaussian sums, $(2.39)$ follows immediately 
from $(2.30)$.
\medskip
\noindent
{\bf 2.9} In this section we shall look into the relation between
the Large Sieve and Selberg's Sieve, in a perspective 
different from the above;
we shall show that the multiplicative Large Sieve can be
amalgamated with Selberg's Sieve. As the discussion of the
previous section suggests, such an extension of the multiplicative
Large Sieve has consequences in the theory of the distribution of
primes in arithmetic progressions$\,{}^{41)}$. This aspect should not be
unexpected, especially if it is taken into account that an origin of
Selberg's Sieve can be traced back to $(2.37)$. In fact, an initial version
of Selberg's  procedure developed in $(2.13)$--$(2.19)$ could be found in
his argument to compute the minimum value of the expression $(2.37)$
under the side condition $a_1=1$$\,{}^{42)}$; 
that is, the extremal values of
$a_n$ are found in much the same way as those of $\lambda(d)$.
\par
Returning to $(2.19)$, the optimal $\lambda$ is written as
$$
\sum_{n\in\Omega(d)}\lambda(d)={1\over G(z,\Omega)}
\sum_{q< z}\mu(q)^2H(q,\Omega)\Psi_q(n,\Omega),\quad
\Psi_q(n,\Omega)=\prod_{\scr{n\in\Omega(p)}\atop\scr{p|q}}
\left({-1\over H(p,\Omega)}\right).\leqno(2.40)
$$
We compare this with $(2.30)$, and ponder upon
the norm of the linear operator $\left(\psi_q(n,\Omega)\right)$, with
$$
\psi_q(n,\Omega)=\mu(q)
\sqrt{H(q,\Omega)}\Psi_q(n,\Omega).\leqno(2.41)
$$
In this way, we find that for any complex vectors 
$\{a_n\}$, $\{b_q\}$$\,{}^{43)}$ 
$$
\leqalignno{
\sum_{q< z}\left|\sum_{M\le n< M+N}
a_n\psi_q(n,\Omega)\right|^2&\le(N-1+z^2)
\sum_{M\le n< M+N}|a_n|^2,&(2.42)\cr
\sum_{M\le n< M+N}\left|\sum_{q< z}
b_q\psi_q(n,\Omega)\right|^2
&\le(N-1+z^2)\sum_{q< z}|b_q|^2.&(2.43)
}
$$
Setting $a_n=\varpi(n)$ in $(2.42)$ we get $(2.10)$ again; and $(2.43)$
implies $(2.10)$ as well, for it contains $(2.28)$. More generally,
the operator $\left(\chi(n)(k/\varphi(k))^{1/2}
\psi_q(n,\Omega)\right)$ could be viewed in the same way, where $\chi$
is primitive ${\bmod\,k}$, and $(k,q)=1$, $kq< z$. That is, Selberg's Sieve
and the multiplicative Large Sieve could be hybridized, which yields 
interesting refinements of $(2.10)$.
\par
This does not exhaust the flexibility hidden in Selberg's Sieve. An aspect
in which Selberg's Sieve supersedes Linnik's is in that the class of
sequences to which the former is applicable is definitely wider than
that with the latter. For instance, with a given arithmetic
function $f$ one may consider the quadratic form$\,{}^{44)}$  
$$
\sum_{n=1}^\infty
f(n)\left(\sum_{\scr{d|n}\atop\scr{d<z}}
\lambda(d)\right)^2\leqno(2.44)
$$
on the side condition $\lambda(1)=1$. The optimal $\lambda$
thus obtained yields an analogue of the above $\psi_q(n,\Omega)$. It
can be used to extend the multiplicative Large Sieve, which in turn
has an important application; that is, a highly simplified proof
of the Least Prime Number Theorem. In the previous section we stressed
the r\^ole played by the Brun--Titchmarsh theorem in Linnik's proof
of his Least Prime Number Theorem. There the Sieve Method was somewhat
hidden. In the new proof, the Sieve Method emerges as the protagonist, 
and leads the whole story$\,{}^{45)}$.
\vskip 0.7cm
\noindent
{\bf Chapter 3. Rosser's Sieve}
\medskip
\noindent 
{\bf 3.1} In the present chapter we shall return to the circle of
Brun's ideas. Being combinatorial in its nature, Brun's
Sieve demands  efforts to comprehend. On the other hand,
Selberg's Sieve is simple and powerful; also Linnik's Sieve gave rise
to the principle of the Large Sieve, which brought a tremendous impact
to the development of the theory of the distribution of 
primes. Perhaps because of this, 
it took considerably long time
for Brun's theory to be appreciated and shared by many. In fact, it was 30
years later since his work [12] when Rosser (ca.\ 1950) opened a way
leading to the complete settlement of 
the linear sieve$\,{}^{46)}$. Namely,
he discovered a choice of sieve weights on 
the general condition introduced
in Section 1.7 (with $\kappa=1$), which gives best possible main
terms in both the upper and lower bounds. Moreover, the construction
of his sieve weights is relatively simple. In what follows we shall
describe the salient points of Rosser's Sieve, especially his Linear
Sieve. We shall employ symbols and definitions introduced in Chapter
$1$, without mention. We stress that we shall start with 
$\c{S}(\c{A},z_0)$ $(2\le z_0<z)$ instead of $\c{A}$. The reason why
we first sift $\c{A}$ with primes less than a certain $z_0$ will become
apparent in the course of discussion.
\medskip
\noindent
{\bf 3.2} To get a lower bound of the size of a subset in a
finite set, it suffices to have an upper bound of its complementary
subset. To wit, lower bounds could result from upper bounds. This
trivial principle was first exploited effectively by 
A.A. Buchstab [13], in the context of the Sieve Method. More explicitly,
his idea relies on the following identity:
Classifying the elements of ${\cal
S}({\cal A},z_0)\setminus{\cal S}({\cal A},z)$ according to their least
prime factors, we get
$$
|{\cal S}({\cal A},z)|=|{\cal S}({\cal A},z_0)|
-\sum_{z_0\le p<z}|{\cal S}({\cal A}_p,p)|.
\leqno(3.1)
$$
This logical identity is named after Buchstab. If we put $z_0=2$ and
iterate the identity infinitely, we get $(1.2)$. It is, however,
useless in general, and thus Brun introduced a system of
restricting the participating divisors of $P(z)$.
\par
Any restriction of the divisors is the same as to attach the wight $0$
or $1$ to each divisor. With this observation in mind, we 
reconsider $(1.2)$. Thus, let $\eta$ be an arbitrary function with
$\eta(1)=1$, and rewrite $(3.1)$ as
$$
|{\cal S}({\cal A},z)|=|{\cal S}({\cal A},z_0)|
-\sum_{z_0\le p<z}\eta(p)|{\cal S}({\cal A}_p,p)|
-\sum_{z_0\le p<z}(1-\eta(p))|{\cal S}({\cal A}_p,p)|\leqno(3.2)
$$
This is the case with $\ell=1$ of the identity
$$
\eqalign{
\left|{\cal S}({\cal A},z)\right|&=
\sum_{\scr{d|P(z_0,z)}\atop\scr{\nu(d)<\ell}}
\mu(d)\rho(d)|{\cal S}({\cal A}_d,z_0)|+(-1)^\ell
\sum_{\scr{d|P(z_0,z)}\atop
\scr{\nu(d)=\ell}}\rho(d)|{\cal S}({\cal A}_d,p(d))|\cr
&+\sum_{\scr{d|P(z_0,z)}\atop\scr{\nu(d)\le\ell}}\mu(d)\sigma(d)
|{\cal S}({\cal A}_d,p(d))|.
}\leqno(3.3)
$$
Here $P(z_0,z)=P(z)/P(z_0)$, $\rho(1)=1$, $\sigma(1)=0$, and
for $d=p_1p_2\cdots p_l$ $(p_1>p_2>\cdots>p_l)$
$$
\rho(d)=\eta(p_1)\eta(p_1p_2)\cdots\eta(p_1p_2\cdots p_l),\quad
\sigma(d)=\rho(d/p(d))-\rho(d)\quad (p(d)=p_l).
\leqno(3.4)
$$
To prove $(3.3)$, we apply to $(3.2)$ the replacements 
${\cal A}\mapsto{\cal A}_d$,
$z\mapsto p(d)$, $\eta(p)\mapsto\eta(dp)$, and insert the
result into $(3.3)$; then we get $\ell\mapsto\ell+1$. Hence,
setting $\ell>\pi(z)$ in $(3.3)$, we obtain
$$
\left|{\cal S}({\cal A},z)\right|=\sum_{d|P(z_0,z)}
\mu(d)\rho(d)|{\cal S}({\cal A}_d,z_0)|
+\sum_{d|P(z_0,z)}\mu(d)\sigma(d)
|{\cal S}({\cal A}_d,p(d))|,\leqno(3.5)
$$
which is an extension or rather a refinement of $(1.2)$.
\medskip
\noindent
{\bf 3.3} For the sake of simplicity, we impose the restriction
$0\le\eta(d)\le1$ for any $d|P(z)$; thus,
$0\le\rho(d)\le1$, $0\le\sigma(d)\le1$. With this, we shall
try to derive from $(3.5)$ as sharp as possible upper and lower bounds 
of $\left|{\cal S}({\cal A},z)\right|$. First, we observe trivially
$$
(-1)^r\left\{|{\cal S}({\cal A},z)|-
\sum_{d|P(z_0,z)}
\mu(d)\rho(d)|{\cal S}({\cal A}_d,z_0)|\right\}\le
\sum_{\scr{d|P(z_0,z)}\atop\scr{\nu(d)\equiv r\smymod 2}}\sigma(d)
|{\cal S}({\cal A}_d,p(d))|.\leqno(3.6)
$$
There is a way to have the equality here; that is,
$\sigma(d)=0$ $(\nu(d)\equiv r+1\mymod2)$. Thus we set
$$
\eta(d)=1\quad(\nu(d)\equiv r+1\mymod2);\leqno(3.7)
$$ 
and such an $\eta$ is denoted as $\eta_r$, and correspondingly
we define $\rho_r$ and $\sigma_r$. Then, we have
$$
|{\cal S}({\cal A},z)|=
\sum_{d|P(z_0,z)}
\mu(d)\rho_r(d)|{\cal S}({\cal A}_d,z_0)|+(-1)^r
\sum_{d|P(z_0,z)}\sigma_r(d)
|{\cal S}({\cal A}_d,p(d))|.\leqno(3.8)
$$
Discarding the second sum, we get
$$
(-1)^r\left\{|{\cal S}({\cal A},z)|-
\sum_{d|P(z_0,z)}
\mu(d)\rho_r(d)|{\cal S}({\cal A}_d,z_0)|\right\}\ge0.\leqno(3.9)
$$
\par
Also we have, corresponding to $(3.8)$,
$$
V(z,\omega)=V(z_0,\omega)V_0(z,\omega;\rho_r)
+(-1)^r\sum_{d|P(z_0,z)}\sigma_r(d){\omega(d)\over d}V(p(d),\omega),
\leqno(3.10)
$$
where 
$$
V_0(z,\omega;\rho_r)=\sum_{d|P(z_0,z)}\mu(d)\rho_r(d)
{\omega(d)\over d}.\leqno(3.11)
$$
In fact, expanding out the product $V(z,\omega)/V(z_0,\omega)$
and classifying the resulting terms in Buchstab's fashion, we
obtain an identity analogous to $(3.1)$. The rest of discussion
is the same as above. In passing, we note that
$V_0(z,\omega;\rho_r)=V(z,\omega;\rho_r)$, provided $z_0=2$.
\medskip
\noindent
{\bf 3.4} In deriving $(3.9)$ from $(3.8)$ we brought in a certain 
inaccuracy, which should certainly be evaded 
as much as possible$\,{}^{47)}$. For this sake, we note the trivial 
but crucial fact that $|{\cal S}({\cal A}, z)|$ is 
a non-increasing function of $z$. Thus the negligence of 
${\cal S}({\cal A}_d,p(d))$ with $p(d)$ which is small for $\c{A}_d$
causes most likely a relatively large loss. To avoid this we should
better set $\sigma_r(d)=0$ for such $d$. One of the most fruitful
device to make explicit the smallness of $p(d)$ for $\c{A}_d$ is
to introduce two parameters $\beta> 1$ and $D>0$, and to define $p(d)$ to
be small for $\c{A}_d$ if $p(d)<(D/d)^{1/\beta}$. Behind this criterion
is the concept of the Sieving Limit, but at this moment
there is no particular necessity to know the details$\,{}^{48)}$.
\par
Hence, in addition to $(3.7)$, we impose
$$
\eta_r(d)=\cases{1& $p(d)^\beta d<D$,\cr
0&$p(d)^\beta d\ge D$,
}\quad (\nu(d)\equiv r\mymod 2).\leqno (3.12)
$$
Then, $\rho_r$ and $\sigma_r$ are, respectively, 
the characteristic functions of the sets$\,{}^{49)}$
$$
\leqalignno{
{\cal D}(\rho_r)&=\{1\}\cup\left\{d:\, p_1p_2\cdots p_{2k+r-1}
p_{2k+r}^{\beta+1}<D,\, 1\le 2k+r\le l\right\},&(3.13)\cr
{\cal D}(\sigma_r)&=\left\{d:\, \rho_r(d/p(d))=1,\,
p_1p_2\cdots p_{l-1}
p_l^{\beta+1}\ge D,\, l\equiv r\mymod 2\right\},&(3.14)
}
$$
with $r=0,\,1$, and $d$ being the same as in $(3.4)$.
\par
With the sieve weight $\rho_r(d)$ thus constructed, the formula $(1.18)$
is called Rosser's Sieve. In fact, the validity of $(1.18)$ is immediate
in view of $(3.9)$ with $z_0=2$. Moreover, because of the supposition 
$\beta>1$ the level condition $(1.21)$ is fulfilled
with the present $D$. It should be noted
that as $D$ and $\beta$ are taken larger and smaller, respectively, 
the set ${\cal D}(\sigma_r)$ becomes narrower; that is, the loss caused
at the step $(3.9)$ should decrease.
\medskip
\noindent
{\bf 3.5} As an application of Rosser's Sieve, we shall prove
Brun's theorem $(1.22)$ briefly$\,{}^{50)}$. In this section, we work with
$z_0=2$. We note first that if $z^2\le D$, then we have 
$$
{1\over2}\left({\beta-1\over\beta+1}
\right)^{\nu(d)/2}\log D<\log (D/d)\quad 
\left(d\in{\cal D}(\rho_r)\right).\leqno(3.15)
$$
In fact, if $r=1$, $\nu(d)=2\ell$, $\rho_1(d)=1$, then we have
$p_{2j+2}<p_{2j+1}<\left(D/(p_1p_2\cdots p_{2j})\right)^{1/(\beta+1)}$ 
$(0\le j\le\ell-1)$. Thus, $((\beta-1)/(\beta+1))\log (D/(p_1p_2\cdots
p_{2j}))<\log (D/(p_1p_2\cdots p_{2j+2}))$, which gives $(3.15)$. Other
cases are analogous. Next, in the second sum of $(3.10)$, the terms are
classified according to the values of $\nu(d)$, and taking 
$(1.17)$ into account we see that the sum is

$$
\ll V(z,\omega)\sum_{l=\ell}^\infty{1\over l!}
\left({\log z\over\log q}\right)^\kappa
\left(\sum_{q\le p<z}{\omega(p)\over p}\right)^l,\leqno(3.16)
$$
where $q=\min_{\nu(d)=l}p(d)$, $\ell=\min\nu(d)$ with $d|P(z)$ and
$\sigma_r(d)=1$. By the definition $(3.4)$, $\rho_r(d)=0$, and thus
$p(d)^\beta d\ge D$, which gives $\ell\ge \tau-\beta$ because 
$z^\tau=D$. On the other hand, we have $\rho_r(d/p(d))=1$, and by $(3.15)$ 
$$
{1\over2}\left({\beta-1\over\beta+1}
\right)^{(\nu(d)-1)/2}\log D<\log(D/d)+\log p(d)\le
(\beta+1)\log p(d),\leqno(3.17)
$$
which gives a lower bound for $q$. Inserting these assertions on $\ell$
and $q$ into $(3.16)$, and setting $\beta=\tau/3$, we reach $(1.22)$ after
some elementary estimation.
\medskip
\noindent
{\bf 3.6} As a matter of fact, it is known that if $\kappa>1$, then
the upper bound  via Rosser's Sieve is inferior to that via
Selberg's Sieve$\,{}^{51)}$. Nevertheless, if $\kappa=1$, then 
Rosser's Sieve yields optimal upper and lower bounds as has been
stressed above. Since the linear sieve problems include great
conjectures, Rosser's construction of his Linear Sieve is
extremely important$\,{}^{52)}$.
\smallskip
We first show his assertion: With $\beta=2$, we fix Rosser's sieve weights 
$\rho_r(d)$; and let the functions $\phi_r(\tau)$,
with $\phi_{r_1}\equiv\phi_{r_2}$ $(r_1\equiv r_2\mymod 2)$, satisfy the
difference-differential equation$\,{}^{53)}$
$$
\leqalignno{
&{d\over d\tau}\left(\tau\phi_r(\tau)\right)
=\phi_{r+1}(\tau-1)\quad
(\tau\ge2),&(3.18)\cr
&\tau\phi_1(\tau)=2e^{c_E},\quad \phi_0(\tau)
=0\quad (0<\tau\le2).
&(3.19)
}
$$
Then we have, for
$z=D^{1/\tau}$, $z_0=\exp((\log D)/(\log\log D)^2)$,
$$
V(z_0,\omega)V_0(z,\omega;\rho_r)
=(1+o(1))\phi_r(\tau)V(z,\omega)
\quad (\kappa=1).\leqno(3.20)
$$
\smallskip
With this, we apply Brun's Sieve $(1.22)$ to the term 
$|{\cal S}({\cal A}_d, z_0)|$ appearing in $(3.9)$, and find that
$$
(-1)^r\left\{|{\cal S}({\cal A},z)|
-(1+o(1))\phi_r(\tau)V(z,\omega)X
\right\}\ge (-1)^r\sum_{\scr{d<D}
\atop\scr{d|P(z)}}\mu(d)\delta_r(d)R_d,\leqno(3.21)
$$
in which  $D=z^\tau$, and $\delta_r$ is a certain characteristic 
function$\,{}^{54)}$. This is called Rosser's Linear Sieve. 
\smallskip
The main steps of the proof of $(3.20)$ are as follows: With $\beta>1$, which
is to be fixed later, we construct Rosser's sieve weights $\rho_r(d)$,
and put 
$$
V(z,\omega)K(z,\rho_r)=\max\left\{0, V(z_0,\omega)
V_0(z,\omega;\rho_r)\right\}.\leqno(3.22)
$$
Then we assume that there exist continuous functions
$k_r$ such that $K(z,\rho_r)=(1+o(1))k_r(\tau)$; we have
$0\le k_0(\tau)\le1\le k_1(\tau)$ by $(3.10)$. Also, 
we assume$\,{}^{55)}$
that $\beta=\inf\{\tau:\,k_0(\tau)>0\}$. 
On noting the definitions $(3.11)$ and $(3.13)$, we have
$$
\eqalign{
V_0(z,\omega;\rho_1)
=V_0(z_0,\omega;\rho_1)-\sum_{z_0\le p<z}{\omega(p)\over p}
V_0(p,\omega;\rho_0^*).
}\leqno(3.23)
$$
Here the level of the Rosser sieve weight $\rho_0^*$ is 
equal to $D/p$; that is, $p_1=p$ in $(3.13)$. If $\rho_1(p)=1$, then 
$p<D^{1/(\beta+1)}$. Thus, if $z^{\beta+1}\ge D$, then 
$V_0(z,\omega;\rho_1)=V_0(D^{1/(\beta+1)},\omega;\rho_1)$.
In view of $(1.17)$ ($\kappa=1$), we have 
$\tau k_1(\tau)=(\beta+1)k_1(\beta+1)$ $(\tau\le\beta+1)$.
If $\tau>\beta+1$, then $\log(D/p)/\log p>\beta$; thus,
by our assumption on $\beta$, we may write
$V(z_0,\omega)V_0(p,\omega;\rho_0^*)=V(p,\omega)
K(p,\rho_0^*)$. Hence, we find that
$$
V(z,\omega)k_1(\tau)
=V(z_0,\omega)k_1(\tau_0)-(1+o(1))\sum_{z_0\le p<z}
{\omega(p)\over p}V(p,\omega)k_0
(\xi_p),\leqno(3.24)
$$
with $z_0^{\tau_0}=D$ and $\xi_p=(\log D)/(\log p)-1$.
We apply $(1.17)$ to $(3.24)$, and express the result 
in terms of a Stieltjes integral. We are led to the
integral equation
$$
\tau_0k_1(\tau_0)-\tau k_1(\tau)=
\int_\tau^{\tau_0}k_0(\xi-1)d\xi.\leqno(3.25)
$$
This ends the discussion on the case $r=1$. The other case
could be treated analogously. The equation that $k_0$
should satisfy is $(3.18)$ if $\tau\ge\beta$; and if
$\tau<\beta$, then it is $(3.19)$ but with the constant
$2e^{c_E}$ being replaced by $(\beta+1)k_1(\beta+1)$.
From this, it follows readily that
$$
k_r(\tau)=1+(-1)^r{(\beta-2)\over\tau^2}(\beta+1)
\phi_1(\beta+1)\left(1+O\left({1\over\tau}\right)\right)
+O\left({1\over \Gamma(\tau)}\right).\leqno(3.26)
$$
In view of Brun's theorem $(1.22)$, we find the optimal
value of $\beta$; that is, we should set $\beta=2$. 
Then, $3k_1(3)=2e^{c_E}$ follows. In this way we reach
$(3.18)$--$(3.19)$.
\par
Having this, we set $k_r=\phi_r$, and go from $(3.25)$ back to
$(3.24)$; then we see that $(3.24)$ holds in fact with $k_1=\phi_r$ and 
$k_0=\phi_{r+1}$. By the definition $(3.19)$ one may multiply
each summands on the right by the factor $\rho_r(p)$. The
identity thus obtained can be iterated in much the same way
as in Section 3.2. We get
$$
V(z,\omega)\phi_r(\tau)
=(1+o(1))V(z_0,\omega)\sum_{d|P(z_0,z)}
\mu(d)\rho_r(d){\omega(d)\over d}
\phi_{r+\nu(d)}\left({\log(D/d)\over\log z_0}
\right).\leqno(3.27)
$$
Comparing this with the definition $(3.11)$ and on noting
$\phi_r(\tau)=1+O(1/\Gamma(\tau))$ ($(3.26)$, $\beta=2$), 
our problem is reduced to the estimation of the expression
$$
\sum_{d|P(z_0,z)}\rho_r(d)
{\omega(d)\over d}
\exp\left(-{\log(D/d)\over\log z_0}
\right)\leqno(3.28)
$$
We skip the details, but this can be seen to be negligible,
which ends the proof of $(3.20)$.
\medskip
\noindent
{\bf 3.7}
The extremal situation$\,{}^{56)}$ that implies that 
Rosser's Linear Sieve is optimal is given by
$$
{\cal B}^{(r)}=\left\{n<x :\, \hbox{[the total number
of prime divisors of $n$] $\equiv 
r\bmod\, 2$}\right\}.\leqno(3.29)
$$
We have $X=x/2$ and $\omega\equiv1$.
Rosser's Sieve $(\beta=2,\,D=x)$ gives
$$
|{\cal S}({\cal B}^{(r)},z)|=\sum_{d|P(z)}\mu(d)
\rho_r(d)|{\cal B}_d^{(r)}|.\leqno(3.30)
$$
Hence, no loss is caused at $(3.9)$ with the present specialization. 
The argument of the previous section could be repeated, and we get
$$
|{\cal S}({\cal B}^{(r)},z)|=(1+o(1))
{x\over2}\phi_r\left({\log x\over\log z}\right)\cdot V(z,1).\leqno(3.31)
$$
Namely, the upper and lower bound implied by Rosser's Linear Sieve
are in fact attainable.
\medskip
\noindent
{\bf 3.8} As an application of the above, we exhibit 
J.-R.\ Chen's theorem [15]$\,{}^{57)}$: For any sufficiently large even
integer $N$, we have$\,{}^{58)}$
$$
\left|\left\{p:\, N=p+P_2\right\}\right|> C_0
\prod_{p>2}\left(1-{1\over(p-1)^2}\right)
\prod_{\scr{p|N}\atop\scr{p>2}}
\left({p-1\over p-2}\right)
{N\over (\log N)^2},\leqno(3.32)
$$
with an absolute constant $C_0$. Here $P_2$ denotes an
integer which has two prime factors at most. With no doubt,
this famous assertion is at the pinnacle of the entire modern theory of
the Sieve Method.
\smallskip
Chen's plan of the proof is relatively simple. We first 
pick up any integer $n<N$ such that $(n,P(N^{1/10}))=1$, and
consider the value of the expression$\,{}^{59)}$ 
$$
W(n)=1-{1\over2}
\sum_{\scr{p_1|n}\atop\scr{N^{1/10}
\le p_1<N^{1/3}}}1-{1\over2}
\sum_{\scr{p_1|n}\atop\scr{N^{1/10}\le p_1<N^{1/3}}}
\sum_{\scr{n=p_1p_2p_3}\atop
\scr{N^{1/3}\le p_2<(N/p_1)^{1/2}}}1.\leqno(3.33)
$$
We find readily that if $W(n)>0$ or $W(n)\ge {1\over2}$, then $n$ is a
$P_2$.
Thus, with ${\cal A}=\{N-p:\,  p<N\}$, we have that
$$
\eqalign{
&\left|\left\{p:\, N=p+P_2\right\}\right|\ge 
|{\cal S}({\cal A}, N^{1/10})|-{1\over2}\sum_{N^{1/10}\le p_1
<N^{1/3}}|{\cal S}({\cal A}_{p_1}, N^{1/10})|\cr
&-{1\over2}\left|\left\{p<N:\, N=p+p_1p_2p_3,
N^{1/10}\le p_1<N^{1/3}\le p_2<(N/p_1)^{1/2}
\right\}\right|.
}\leqno(3.34)
$$
We apply $(3.21)$ $(r=0)$ to the first term on the right,
and  $(3.21)$ $(r=1)$ to the second. On the other hand, we
replace the third term by 
${1\over2}|{\cal S}({\cal A}^*, N^{1/2}(\log N)^{-A})|$
with $A$ being sufficiently large, and apply $(3.21)$ $(r=1)$ again. 
Here ${\cal A}^*=\{N-p_1p_2p_3:\, p_1p_2p_3<N\}$
with $p_1,p_2$ as above. To bound
the remainder terms that arise from the first two applications of Rosser's
Linear Sieve, we employ Bombieri--Vinogradov's Mean Prime Number Theorem
$(2.38)$. To deal with the remainder term caused by the third application, 
we employ an extension$\,{}^{60)}$ of $(2.38)$ to the sequence
$\{p_1p_2p_3<N\}$ with $p_1,p_2$ as before. 
What remains is to compute the main terms, which is, however, 
rudimentary.
\vskip 0.7cm
\noindent
{\bf Chapter 4. The Remainder Term}
\medskip
\noindent
{\bf 4.1} From what we have described so far, one may infer the
reach of the modern Sieve Method. To continue the story,
we should now leave the discussion of the main terms for
the estimation of the remainder terms$\,{}^{61)}$. 
There must converge the true
essences of Analytic Number Theory, as the proof of Chen's
theorem $(3.32)$ illustrates dramatically. We have, however,
too many relevant fields, subjects, and technicalities to mention. Thus we
would rather single out the idea that is probably the most
fundamental in the theory of the remainder terms, especially of the Linear
Sieve. A culmination in this context is due to H. Iwaniec [34], and to
describe it we need to tell a brief history.
\medskip
\noindent
{\bf 4.2} The development started with the discovery
that Selberg's sieve weights could intervene in the control
of the remainder term, in a highly non-trivial way;
the serendipity occurred to us [49].  This is in fact a
surprising fact,  because those weights had been 
constructed solely with the
aim to attain the best possible main term
while the remainder term had been utterly disregarded. That is,
in the structure of the sieve weights thus defined is hidden a
mechanism that could induce massive cancellations among the summands
in the second sum of $(1.20)$, $r=1$. 
A little later the same occurred to Chen [16] 
with Rosser's Linear Sieve.
\medskip
\noindent
{\bf 4.3}  
Let us dwell a little on their findings. That is about the level
of sieve weights, and we have to return to the concept itself. At the
bottom is the too natural prerequisite that any
main term be superior in magnitude to the corresponding remainder term. 
However, in the Sieve Method this triviality had never been achieved in any
effective manner until Brun introduced the cutoff argument into
Eratosthenes' Sieve, and brought about a revolutionary change. From there
stemmed the concept of the level of sieve weights, as introduced  at
$(1.21)$, though it left for long its trail only in the  primitive 
$$
|R({\cal A},z;\rho_r)|\le\sum_{\scr{d<D}\atop\scr{d|P(z)}}
|\rho_r(d)||R_d|.
\leqno(4.1)
$$
Having the sieve weights
that yield optimal main terms, the
focus of attention is naturally on the very basic issue:  to get larger
levels.  That is essentially the unique way to make the inequality
$(1.19)$ sharper as the reasoning at the end of
Section 3.4 suggests, despite it is meant only for Rosser's Sieve.
Namely, to go beyond $(4.1)$  the inner structure of
the sequence $\{\mu(d)\rho_r(d)R_d\}$  has to be exploited so that the
cancellation among the members be detected; and the size of $D$, the
level of $\rho_r$, could  be taken larger than a priori. 
\par 
This was done for the first time in [49]. Thus, let  $\Omega(p)=\{0\}$,
say, and write $\mu(d)\rho_1(d)=\sum_{[d_1,d_2]=d}
\lambda(d_1)\lambda(d_2)$ with $\lambda(d)=0$ for $d\ge z$. 
We assume that $\lambda$ is chosen optimally 
so that $|\lambda(d)|\le 1$  as $(2.20)$ shows. 
Then Selberg's Sieve or $(2.11)$ implies that $(4.1)$ ($r=1$) could be
replaced by the expression
$$
|R({\cal A},z;\rho_1)|
\le\sup_{a,\,b}\left|\sum_{m<z}\,\sum_{n<z}a_mb_n
R_{[m,n]}\right|,
\leqno(4.2)
$$
where $a=\{a_m\}$, $b=\{b_n\}$ are arbitrary vectors such
that $|a_m|, |b_n|\le1$. This is yet trivial; but its implication
is striking. 
For instance, applying $(4.2)$ to the sequence
${\cal A}=\{n<x:\, n\equiv l\mymod k\}$ $((k,l)=1)$, we
obtain an improvement upon $(2.36)$:
$$
\pi(x;k,l)\le2(1+o(1)){x\over\varphi(k)\log(x/\sqrt{k})}
\quad (k\le x^{6/17}).\leqno(4.3)
$$
That is, $(4.2)$ allows us to utilize the level
$D=(x/\sqrt{k})(\log x)^{-2}$ in place of the trivial
$D=(x/k)(\log x)^{-2}$ which is involved in $(2.33)$. The cause
of this is to have had a bilinear form in $(4.2)$;
that is, to have read the structure of those sieve 
weights as such.
\par
On the other hand Chen [16] exploited the Buchstab identity,
in the case of the Linear Sieve. In $(3.1)$ we put $L^J=z/z_0$ with
an integer $J$, divide the sum over the primes according to the  
covering $[z_0,z)=\cup_{j\le J}I$, $I=[z_0L^{j-1}, z_0L^j)$, and
apply Rosser's Sieve to each ${\cal S}({\cal A}_p,p)$, $p\in I$,
with the level $D/(z_0L^j)$. Provided $J$ is chosen appropriately,
the pair of the main terms remains the same asymptotically, because of
$(3.18)$--$(3.19)$; but the remainder term is bounded by
$$
\sup_{K<z}\sup_{a,\,b}\left|\sum_{K\le p<KL}\,
\sum_{n<D/K}\mu(pn)a_pb_nR_{pn}\right|,
\leqno(4.4)
$$
which is to be compared with $(4.2)$. With this bilinear form, 
Chen could detect
the cancellation inside the remainder term. The effect is well
exhibited in his own application that yielded the assertion$\,{}^{62)}$ 
$P_2\in[x-\sqrt{x},\, x)$ for any sufficiently large $x$. 
In view of the relation between the Riemann Hypothesis 
and the existence of primes in short intervals,  
this is indeed remarkable.
\medskip
\noindent
{\bf 4.4} Now, after the two precursors$\,{}^{63)}$, Iwaniec [34] 
made a true incision into the  subject.
Superseding $(4.2)$ and $(4.4)$, he bounded the remainder
term in $(3.21)$ by
the expression 
$$
 (\log z)\cdot
\sup_{a,\,b}\left|\sum_{m<M}\,\sum_{n<N}\mu(mn)a_mb_n
R_{mn}\right|,\leqno(4.5)
$$
where $M,N$ are arbitrary except for $MN=D$. This is called Iwaniec's
bilinear form for the remainder term in the Linear Sieve. The basis of his
idea is fairy simple; that simplicity is shared by the above two
ideas similarly. Thus, in the process to reach Rosser's sieve weights
$(3.13)$, $\beta=2$, those primes participating the sieve are
classified as Chen did. The function $\eta$ is first
defined over the family of all set theoretic products
of intervals $I$; and it is redefined as a function over
integers, in an obvious manner 
according to their prime decompositions.
With this, we proceed in much the same way as we did 
in Sections 3.2--3.4
and 3.6. Then a smoothed version of $(3.21)$ emerges. 
What remains is solely
Iwaniec's penetrating observation on the remainder term 
thus obtained$\,{}^{64)}$.
\par
That the parameters $M,N$ are independent is a real merit
in Iwaniec's Linear Sieve, because of which $(4.5)$ has given rise to many
remarkable consequences. One of the best applications is done by Iwaniec and
M. Jutila [36], a landmark among the works on the existence of
primes in short intervals$\,{}^{65)}$.
\vskip 0.7cm
\noindent
{\bf Conclusion} 
\medskip
\noindent
What $(4.5)$ for instance suggests is the importance of
the circle of methods, which are represented by Linnik's Dispersion
Method. The origin could be found in the Weyl--van der Corput method
dealing with trigonometrical sums, which is a device closely
related to subconvexity bounds of the Riemann zeta and analogous
functions, though a far cry from the Lindel\"of  Hypothesis. 
\par In those 
methods, especially in Linnik's, often 
Kloosterman sums play a fundamental r\^ole. This was the very reason 
why Linnik [44] envisaged the cancellation among
the sums, which is perhaps hard to detect if one sticks to
algebraic means only. Selberg [71] opened a way, and V.N. Kuznetsov
[40] made a remarkable contribution to realize a part of 
Linnik's dream. Together with an independent research by
R.W. Bruggeman [8], the work [40] brought a new era in Analytic
Number Theory. That was began by Iwaniec [32], when he combined
their works with the additive Large Sieve and created the spectral
Large Sieve. On that basis, Bombieri, J.B. Friedlander and Iwaniec [7] 
achieved a genuine improvement upon $(2.38)$.
\par
On the other hand, the appearance of Kloosterman sums in the discussion 
of non-diagonal parts arising in the applications of the Dispersion Method 
or alike must be a reflection of the fact that we are actually working
on a certain group structure$\,{}^{66)}$, that is, we are looking at the 
remainder term in sieves via automorphic and harmonic mechanisms on
${\rm GL}(2)$. Bearing Brun's torch we have come to a far country.
\bigskip
\noindent
{\bf Addendum} (May 14, 2005)
\medskip
\noindent
Very recently there was a fantastic development in the study of gaps
between primes:  In their unpublished preprint 
`{\it Small gaps between primes.\ II} ({\it preliminary\/})' 
(February 8, 2005),
D.A. Goldston, J. Pintz, and C.Y. Y{\i}ld{\i}r{\i}m established, among other
things,
$$
\liminf_{n\to\infty}{p_{n+1}-p_n\over\log p_n}=0,\leqno(\r{A}.1)
$$
and that if $(1.26)$ holds for $Q=x^\theta$ with
a $\theta>{1\over 2}$ then there exists an absolute constant $c(\theta)$
such that
$$
\liminf_{n\to\infty}(p_{n+1}-p_n)\le c(\theta),\leqno(\r{A}.2)
$$
where $p_n$ is the $n$th prime.
For a short, essentially self-contained proof of these facts, see the
preprint arXiv:math.\ NT/0505300 `{\it Small gaps between primes exist\/}' by
Goldston, Motohashi,  Pintz, and Y{\i}ld{\i}r{\i}m. The argument is in the 
framework of Selberg's Sieve and Linnik's Large Sieve; the latter 
is in the sense that $(2.38)$ plays a fundamental r\^ole. 
One might surmise that a proof of Twin Prime Conjecture is not beyond
the reach of today's Analytic Number Theory.
\bigskip
\noindent
{\bf Addendum} (September 20, 2006)
\medskip
\noindent
Because of the improvement [7] of the Mean Prime Number Theorem, the last 
hypothetical assertion on {\it bounded
differences between consecutive primes\/} might appear within the reach 
of the present technology. Until recently there were, however, two main obstacles
that prevented us to make any real incision into the matter: One is the fact
that Selberg's Sieve or rather his sifting procedure, which is highly essential
in the argument of Goldston, Pintz, and Y{\i}ld{\i}r{\i}m, did not seem
to admit any error terms with the flexibility that $(4.5)$ enjoys, although we had
already a partial result [62] that corresponds to the situation 
$M=N$ in $(4.5)$. Another obstacle is that in [7] only 
the distribution of primes in arithmetic progressions with a {\it fixed} residue
class, i.e., $\pi(x;q,a)$ with $a$ fixed, is considered, and to
achieve $(\r{A}.2)$ we need to relax this restriction to a considerable extent. 
\par
With this, the present situation is that although the latter difficulty still
persists, Motohashi and Pintz have
succeeded in suppressing the former in their preprint arXiv:math.\ NT/0602599 `{\it A
smoothed GPY sieve\/}' by extending the argument of [62]. We note, by the way, that
arXiv:math.\ NT/0505300 quoted above has been published in Proc.\ Japan Acad.\ 82A
(2006), 61--65, which can be downloaded freely via the web-page of the academy.
\bigskip
\noindent
{\bf Acknowledgements}
\medskip
\noindent
On this occasion I would like to express my indebtedness to the
late Prof.\ P. Tur\'an for all that he bestowed upon me; to  
Profs.\ M. Jutila and A. Ivi\'c for their unfailing friendship and constant
encouragement; to my colleagues for their generosities; and
to my family for the comfort.
\par
I thank G. Greaves for the numerous textual improvements 
incorporated in this new version.
\vskip 1cm
\noindent
{\bf Notes}
\medskip
\noindent
\item{1)} My impression that has arisen reading occasionally the
history of the ancient Meso\-potamia. As to the remote origin of the concept
of prime numbers and decomposition of integers into prime factors, there
exist a collection of highly striking evidences from Sumer and Akkad; for instance,
see the recent discovery [63] by K. Muroi. Ancient Babylonian mathematicians
tried to make exercises for their students harder or more enjoyable by
embedding not only the familiar $2, 3, 5$ but also
$$
7, 11, 13, 17, 19, 23, 29, 31, 41,
47, 59, 79, 83, 137, 139, 1481
$$ 
into linear and quadratic equations as hidden prime factors of the coefficients.
They were playing with prime numbers, more than 4000 years ago. To me those primes
appear exactly like glittering gems inlaid into the famous treasures from Ur, Nimrud,
and the grave of King Tutankhamen.
Perhaps more than that, because those treasures are perishable but mathematical
equations are never. Certainly this list of primes will be expanded in the near
future as the research should develop further.
\smallskip\noindent
\item{2)} Despite this common attribution, it appears likely 
that the procedure had come from the Orient.
Ancient mathematicians made academic trips as we do today; 
an Ebla tablet tells such an episode 
(G. Leick.\  Mesopotamia.\ Penguin Books 2002, p.\ 68).
Eratosthenes was a poet,
astronomer, and director of the great library of Alexandria (3rd century
BCE); he is said to have been old and unable to enjoy books 
when he died on voluntary starvation. 
\smallskip
\item{3)} There are two basic queries in the Sieve Method. One is the
primality test, and the other the number of 
primes that satisfy a certain set of conditions. To the former,
N. Agrawal et al made a remarkable contribution a few years ago. However,
my interest lies mainly in the latter query, the quantitative
aspect of the distribution of primes. 
To determine a spot to dig where a ruby could certainly be found
is definitely harder than to determine whether a stone is a ruby 
or not. I am well aware that the word `ineffectual' is too drastic, for
Eratosthenes' Sieve could sometimes be a useful tool. 
A typical instance is in I.M.
Vinogradov's proof (1937) of the ternary Goldbach Conjecture (see [64,
Kap.\ VI] and  R. C. Vaughan's newer argument [76]). See also H. Iwaniec
[29].
\smallskip
\item{4)} Prior to Brun [10] was an incomplete trial by J. Merlin [45]. Brun
mentioned this fact in [11] precisely; an indication of his respect for 
the academic axiom, which is often purposely 
ignored nowadays.  Fortunately, Brun's name seems to be well-known. For
instance, at the item `Number Theory' in World Encyclopedia 
(Heibon, Tokyo 1972; Japanese), C. Chevalley writes `the unique result
known about the Twin Prime Conjecture is the work  by Viggo Brun (1919)' (see
the bottom lines of the right column  on p.\ 458, vol.\ 16). This is a
bizarre opinion, however. Scripta manent.
\smallskip
\item{5)} A.M. Legendre's formulation (1808), though the use of the M\"obius
function (1832) is a later tradition.
\smallskip
\item{6)} This asymptotic formula is an elementary result due to F. Mertens
(1874) (see [64, pp.\ 80--81]); $\exp(-c_E)=0.561459483566885\ldots$.
\smallskip
\item{7)} As to Prime Number Theorem, see, e.g., [26] [64]. 
Here the theorem is used in the abused manner 
that the probability of the appearance of a prime at
$x$ is $1/(\log x)$.
\smallskip
\item{8)} Under the condition $y=x^\theta$, $0<\theta<{1\over2}$, 
this appears impossible to prove even on the Riemann Hypothesis.
See Section 2.8.
\smallskip
\item{9)} The second assertion follows from the first, because of 
$x/(\log x)\ll\pi(x)$, an elementary assertion due to 
P.L. \v{C}eby\v{c}ev (1853). See [64, p.\ 19].
\smallskip
\item{10)}Is it an exaggeration to say that it took more than 2000 years
for number theorists to break themselves of cleaving to equalities?
A great king cut the knot with a stroke of his sword.
\smallskip
\item{11)} See e.g., [20, \S3.4] for details.  
\smallskip
\item{12)} This formulation of the sieve dimension might appear 
made abruptly. I have in fact skipped the discussion about
how to define a mean value of $\omega$. See, for instance, 
[20, pp.\ 27--37]. In practice, $(1.17)$ should suffice.
Note that there is no a priori sieve dimension for a given
problem to which the Sieve Method is to be applied. 
It depends largely upon the choice of the sequence
$\c{A}$ into which the problem is embedded. The two subsequent
sections give examples in this context as well.
\smallskip
\item{13)} Be aware that at this stage nothing is assumed about the
remainder term. Obviously it is useless to employ sieve weights
that make the estimation of the remainder term too hard.
Although Brun [12] stated his result including an estimation
of the remainder term, nowadays the main term is first treated,
and then the discussion about the remainder term follows.
A reason for this splitting of the theory into two parts
could be seen from e.g., Section 1.9 below.
\smallskip
\item{14)} It is proved in [12] that both of the conjectures could be resolved
if one is allowed to replace primes by integers with at most
$9$ prime factors. See Section 3.8 below.
\smallskip
\item{15)} About the distribution of primes in arithmetic progressions,
see [64, Kap.\ IV]. As to the distribution of primes and the
Riemann Hypothesis, see [64, p.\ 235]. If an unconditional 
uniform bound for individual $E(x;k,l)$ is sought for, it is hard to
supersede the famous but ineffective Siegel--Walfisz's 
Prime Number Theorem (see [64, p.\ 144]); namely, $Q=(\log x)^B$ with any
fixed $B>0$ is the best presently, and any assertion like $(1.28)$ may
appear hopeless. See Section 2.8, however.
\smallskip
\item{16)} See e.g., [64, Kap.\ VI].
\smallskip
\item{17)} More generally, one may discuss with the system of residues
$\Omega(p^a)\bmod\, {p^a}$, $a=1,2,\ldots$, being given initially. 
See Selberg [74] and also my lecture notes [60, \S 1.1] as well as [54, III].
\smallskip
\item{18)} Without loss of generality, one may suppose that
$0<|\Omega(p)|<p$, as I do within the present chapter. 
In fact, if $|\Omega(p)|=0$, then such a prime does not participate the
sieve problem under consideration; and if $|\Omega(p)|=p$, then
$|{\cal S}({\cal A};\Omega,z)|=0$ as soon as
$z$ becomes larger than $p$. One could restrict the domain of
$\Omega$ instead.
\smallskip
\item{19)} The restriction $q<z$ is made solely 
for the sake of simplicity. In general, $q<Q$, $q|P(z)$,
should be employed,
with $Q$ depending on $z$.
\smallskip
\item{20)} As a matter of fact, Linnik [41] dealt 
with the case in which $q$ are
all primes. His motivation was to investigate statistically I.M.
Vinogradov's conjecture on least quadratic non-residues 
(an account can be found in [6, p.\ 7]). 
Thus the discussion of the present
section is a refined version of Linnik's argument. The beautiful
inequality $(2.10)$ is in fact due to H.L. Montgomery [46, Chap.\ 3]. 
See also [54, II].
\smallskip
\item{21)} Selberg's argument is indeed rich; that will become apparent
in Sections 2.9 and 4.3 below. 
\smallskip
\item{22)}The estimation $(2.20)$ is not optimal; neither the bound for
$|R|$ in $(2.22)$. In the case where $|\Omega(p)|$ is bounded,
these can be acceptable. However, for instance, 
if $|\Omega(p)|\sim cp$ with a constant
$c>0$, then $(2.19)$ yields a bound of $\lambda(d)$ far superior to
$(2.21)$. In passing, we note that Selberg's Sieve in the context
of the present section can be applied to general sequences
upon the supposition $(1.3)$; see e.g., [20, Chap.\ 2] for details.
\smallskip
\item{23)} Probably because of this comparison, it is said often
that Selberg's Sieve is effective only when $|\Omega(p)|$
is relatively small. This is, however, incorrect, as is proved in
the next section. See also [60, \S1.1] [54, II] [54, III].
\smallskip
\item{24)} The duality between Linnik's and Selberg's Sieves seems
to have been observed for the first time by myself at a Tur\'an seminar
in 1970, as explicitly as is given here. Thus the contents
of this section dates back to 1970, although it was published
in [54, II] much later because of an unfortunate circumstance.
\smallskip
\item{25)} Selberg's inequality; obviously an extension
of Bessel's inequality. Proof is easy; see e.g., 
[6, p.\ 14] [46, pp.\ 42--43]. The factor $N-1+\delta^{-1}$ in
$(2.30)$ is due to Selberg, and is best possible, a proof of which
is in [47]. The history of the Large Sieve, starting at Linnik [41] and
reaching the contents of Section 2.6 (early 1970's), 
is highly interesting. Above all, the great leap made 
by Bombieri [4] should be acclaimed. That impact brought
a number of then young people into Analytic Number Theory; some of them 
are still active. In between Bombieri [4] and Selberg's
inequality is G. Hal\'asz [21].
\smallskip
\item{26)} However, one should not forget the fact that
Brun's Sieve is capable of giving rise to lower bounds as well;
for instance, $(1.23)$.
\smallskip
\item{27)} See [60, pp.\ 129--130]. As to the relation between
the distribution of primes in arithmetic progressions and
the Exceptional Zeros or the Siegel Zeros of Dirichlet $L$-functions,
see [64, Kap.\ IV]. 
\smallskip
\item{28)} This is due to Montgomery [46, Chap.\ 4]. There 
arbitrary intervals are in fact considered. More precisely,
he used instead of $(2.30)$ a sharper inequality 
deducible from $(2.29)$. A completely uniform result was
later obtained by Montgomery--Vaughan [48]. Here is an
important note: Combining the discussion of Section 2.5
with $(2.29)$, one concludes that $(2.17)$ is not optimal. 
That is, it is suggested that the appeal to $(2.29)$
should yield a procedure with which one could aim 
a simultaneous optimization for the main and remainder terms
in Selberg's Sieve applied to intervals. Similar observation
is made on [22, p.\ 126]. However, my older argument developed
in Section 2.5 seems to go deeper.
\smallskip
\item{29)} My usage of the word `avoid' in the present context
needs to be explained. It indicates, in a somewhat abused way, 
that there exist situations where one may
reach significant results on the distribution of primes without
appealing to the great hypothesis, though arguments might become
more involved than on the hypothesis, and results 
less sharp. What is
important is that the Twin Prime and the Goldbach Conjectures are contained
in this category of problems. 
\smallskip
\item{30)} To learn under the two great mathematicians,
A. R\'enyi and P. Tur\'an, I arrived Budapest via railway
in the evening of January 31, 1970. Next day,
I saw a black flag over the entrance of Matematikai 
Kutat\'oint\'ezet on Re\'altanoda. The director R\'enyi died aged 49.
\smallskip
\item{31)} The Zero Density Theory. See e.g., [26] [46] [64] for details.
\smallskip
\item{32)} The multiplicative characters $\{n^{it}: t\in\B{R}\}$ could also
be included into the multiplicative Large Sieve. 
Unfortunately, I have to skip basic contributions 
by Hal\'asz [21] and P.X. Gallagher [18]. See [6] [26] [46] for details.
\smallskip
\item{33)} The zero density theorem of the Linnik type. The discussion
given in [64, Kap.\ X \S2] is as hard as Linnik's original, though
some simplifications by K.A. Rodosskii are claimed. 
H. Davenport reviewed Linnik's articles, and commented, `formidable.'
A comparatively simpler proof is that via Tur\'an's Power 
Sum Method [75]; see [6, \S6]. See also [50][60, Chap.\ V]
\smallskip
\item{34)} To understand this fascinating encounter 
between a sieve result and Dirichlet $L$-functions,
I entered into the Sieve Method. The method could resolve a difficulty
that does not appear to be settled with analytic arguments only. 
See [64, pp.\ 346--347]. The argument via the Power Sum Method 
requires as well the Brun--Titchmarsh theorem; see [6, p.\ 50].
\smallskip
\item{35)} See [64, pp.\ 145--146]. Also [52].
\smallskip
\item{36)} R\'enyi neither proceeded nor stated his result like
this. Nevertheless, what he established is essentially the same as
the Mean Prime Number Theorem as is stated here.
\smallskip
\item{37)}The developments prior to May 1964 are detailed in [2];
see [46, Chaps.\ 15--17] for an account of the later progress
up to 1971.
\smallskip
\item{38)} A.I. Vinogradov's assertion is weaker than Bombieri's 
which is embodied in $(2.38)$. In the context of $(1.25)$--$(1.27)$, they are,
however, of the same strength essentially.
\smallskip
\item{39)} The Dispersion Method is a device to introduce perturbations
into certain arithmetic problems. Statistical arguments are then
applied to the perturbed. Thus, for instance, binary additive
problems could sometimes be transformed into
ternary additive problems which are often more tractable. 
\smallskip
\item{40)} Bombieri used his own; but in the context of $(2.38)$
there is no difference. See his lecture notes [6].
\smallskip
\item{41)} The present section is due to Selberg [73] and myself
[50] [52]--[54] [60, \S1.2--\S1.4]. An origin is in [6, Th\'eor\`eme 7A]
(Selberg), which corresponds to the case $\Omega(p)=\{0\}$. 
\smallskip
\item{42)} More correctly, to minimize the main term of the
asymptotic formula for $(2.37)$. See Selberg [66] [67].
\smallskip
\item{43)} A proof is in [60, \S1.2]; see also [54, II] [54, III].
Following Selberg [73],
$\{\Psi_q(n,\Omega)\}$ could be termed pseudo-characters,
which could be regarded as generalizations of the Ramanujan sum.
This relation between these arithmetic functions and 
Selberg's Sieve was found by myself [54, II] [60, Notes(I)]. 
\smallskip
\item{44)} This is just a formal discussion. In practice,
we have to impose realistic conditions to $f$. See e.g., 
[53] [57, Lemma 2] [60, \S1.3].
\smallskip
\item{45)} Due to myself [57]; see also [60, \S6.2]. The two basic 
assertions of Linnik mentioned in Notes 36--38 and Tur\'an's Power Sum
Method are altogether discarded. See also [50][52][53]. These should be
compared with M. Jutila [38] and S. Graham [19]. 
Further, the works [55] [56] [59] are relevant.
\smallskip
\item{46)} Only a few people seem to have had opportunities to look into
Rosser's unpublished manuscript. Other
people, including myself, could see its outline only in Selberg's
scant account [72]. Thus all published works either 
on Rosser's Sieve or on the lower sieve bound, which is implicitly
the main subject of the present section, can be regarded
as contributions made independently of Rosser. With this understanding, 
the first general result on the lower bound 
is in Ankeny--Onishi [1], which
is a combination of Selberg's Sieve and Buchstab's Identity $(3.1)$; see
Note 23 above.  Buchstab [14] is on the line of his [13],
and certainly more sophisticated. The first
published account of  Rosser's Sieve is Iwaniec [33], 
entirely due to himself, a detailed version of which 
is given in [20, Chaps.\ 3--4]. 
As to the Linear Sieve, see Note 52 below.
\smallskip
\item{47)} This section is an excerpt from [60, \S2.1]. See also [58, I].
\smallskip
\item{48)} An easy account of the Sieving Limit is 
in [60, pp.\ 57--60]. For more details see Selberg [72]; and also
Note 61 below. 
\smallskip
\item{49)} These sets are independent of $z_0,\,z$.
\smallskip
\item{50)} This section is due to an observation by
Friedlander--Iwaniec [17]. See also [60, pp.\ 55--59].
The sum over $p$ in $(3.16)$ can be estimated via $(1.17)$.
\smallskip
\item{51)} In this respect, the argument originating in Ankeny--Onishi [1] 
has a definite merit. There exists a detailed discussion in [20, Chap.\ 7]. 
\smallskip
\item{52)} However, the first published account
of the Linear Sieve, i.e., the determination of the optimal upper/lower
main terms is due to Jurkat--Richert [37]. They started with
Selberg's Sieve and performed iteration via $(3.1)$ in much the same
way as Rosser's (i.e., $\beta=2$). On the other hand
Iwaniec [27] is the first published account of the
Linear Sieve \`a la Rosser. Iwaniec [28] dealt with 
the half dimensional sieve, where Rosser's construction with $\beta=1$
yields again the optimal upper/lower main terms.  See also
[20, \S4.5]. By the way, the work of Jurkat and Richert
was applied by J.-R.\ Chen in his famous work [15] on Goldbach's
Conjecture; see Note 60 below. All works quoted here are independent
of Rosser's.
\smallskip
\item{53)} The appearance of this equation might be
somewhat unexpected. It is, however, a consequence of the prerequisite
that the pair of the optimal main terms in the upper and lower bounds
be stable against the iteration via the Buchstab
identity. By the way, $(3.19)$ ($r=1$) coincides with
what Selberg's Sieve implies. Jurkat--Richert [37] starts
with this fact. As to the analogue of $(3.18)$--$(3.19)$ for
general $\kappa$, see [20, \S4.2]. 
\smallskip
\item{54)} Here is the reason why the additional
parameter $z_0$ has been introduced. The discussion in
the rest of this section involves in fact a convergence
issue, though it is not checked there. The r\^ole of $z_0$ is to secure the
convergence. On the other hand, the appearance of the factor
$V(z_0,\omega)$ is  harmless because of Brun's theorem $(1.22)$. In this
context, $(1.22)$ is termed the Fundamental Lemma. The assertion $(3.21)$
results from the combination of Rosser's Sieve and Brun's Sieve 
(or rather the procedure of Section 3.5); the sieve weights of the
former are multiplied by those of the latter, and the new sieve
weights thus obtained are again values of a
characteristic function. As to the condition $d<z^\tau$, it is
in fact the result of taking anew the value of $z$ for a
cosmetic purpose; this is
possible because of the smoothness of $\phi_r$. See [60, Chap.\ III]
for more details.
\smallskip
\item{55)} A heuristic explanation: 
Let $\beta_0$ be the infimum thus defined. Assume that there exists
a summand in the second sum
on the right of $(3.8)$ such that
$d<(D/p(d))^{1/\beta_0}$. Then there exists the possibility that
$|\c{S}(\c{A}_d,p(d))|$ is positive, as this might be detected
by Rosser's Sieve with $\beta=\beta_0$ and the level $D/d$. That is,
$\sigma_r(d)$ has to vanish; and we are led to the 
condition.
\smallskip
\item{56)} Due to Selberg. See [20, \S4.5] for details. 
\smallskip
\item{57)} Most probably, it was early in winter 1966. There was
a colloquium talk by S. Uchiyama at Tokyo University. After his
talk, I had a short discussion with him. He said, 
`An astounding announcement has
been made in China.' `What is that?' 
`A mathematician named Chen Jing-run
has claimed $p+P_2$ for Goldbach's Conjecture.' `Anything about
his method?' `The Mean Prime Number Theorem and two sieve lemmas,
but I have difficulties with one of the latter.' Then he gave me
a copy of the now famous announcement. The dreadful Cultural
Revolution had already been spreading, and Chen would lose seven years
until the publication of the proof.
\smallskip
\item{58)} Chen got $C_0>0.67$. A conjecture of Hardy--Littlewood [23]
states that with $C_0=2$ the right side of $(3.32)$ is asymptotically
equal to $|\{p: N=p+p'\}|$, where $p'$ is also a prime.
\smallskip
\item{59)} The procedure $(3.34)$ is a typical instance of  
weighted sieves, which originates in P. Kuhn [39]. See [20, Chap.\ 5] 
for a general theory.
\smallskip
\item{60)} Chen applied Jurkat--Richert [37] to the first two terms
as mentioned above, and Selberg's Sieve to the third, which is not much
different from the present procedure, as far as $(3.32)$ is concerned.
The move to the sequence ${\cal A}^*$ is now
called the Switching Trick. The extension of Bombieri--Vinogradov's
Mean Prime Number Theorem to $\c{A}^*$ is due to Chen himself.
For a more general extension, see [51] as well as [6, \S22].
\smallskip
\item{61)} Naturally one could take up topics
that are not included in the above category of sieves. 
For instance, there is a sieve method started by Bombieri [5],
which is related to the elementary proof of the Prime
Number Theorem ([64, Kap.\ III \S6]). Recently Friedlander and
Iwaniec (1998) made an important progress on this line.
\smallskip
\item{62)} See [20, pp.\ 257--258].
\smallskip
\item{63)} The seminal nature of the works [16] and [49] could be
stressed with a fair reason. There are explicit mentions
in Iwaniec [31] [34] [35]. 
Certain personal recollections might be
allowed here, as this is presumably the last opportunity for me
to write down some memorable events from my young days: 
The preprint of [49] was finished in early autumn of 1972. There had been a 
belief in the air that
the Brun--Titchmarsh theorem would not be improved beyond 
$(2.36)$ via the Sieve Method. Therefore, I sent copies 
of my works to Chen, Gallagher, 
Halberstam, Hooley, and Richert. Hooley replied me immediately, 
kindly showing how to gain a further improvement. He
had developed a statistical study of the Brun--Titchmarsh theorem.
My preprint seems to have circulated widely, 
with misunderstandings as well;  one day I received a letter from US
informing me that I was rumored to have  proved 
Twin Prime Conjecture. Richert kindly invited me to an MFO Tagung
(1975). After attending the meeting, I went to Budapest to
see my mentor Tur\'an. He would die aged 66
in September next year. He indicated very faintly about 
his illness but no sign of graveness. He wholeheartedly encouraged me 
at the restaurant Astoria after my talk at the institute
on my sieve results including the Brun--Titchmarsh theorem, 
the zero density of the Linnik type, 
and the Deuring--Heilbronn phenomenon. He liked all,
especially the last, even though that made obsolete an important work 
of his already late collaborator S. Knapowski, and thus his Power Sum Method
to a certain extent.
I continued the trip to Warszawa via Krak\'ow to see Iwaniec. 
He and A. Schinzel kindly came to pick up me at
the central station in that early morning of cold December. 
Iwaniec had the ambition to
extend [49] to Rosser's Sieve. At the next MFO Tagun (1977)
he disclosed to me the surprise that he had already got the essentials of his
revolutionary  work [34]. In 1979 my daughter was born. In that summer 
I was at the great Durham Symposium, 
and D. Hejhal kindly gave me a private lecture on
Kuznetsov's work [40], an enormous change of the landscape. There I met also L.-K.
Hua; I told him I wanted to go to Peking to see Chen. Next autumn I could visit
Peking but not Chen because of an unknown reason. In the spring of 1981, I
finished my lecture notes [60] at the Tata IFR. I see still the
magnificent sunset over the Arabian Sea. 
\smallskip
\item{64)} As a matter of fact, we need certain conditions
to have the assertion $(4.5)$ valid. See [58, II] [60, \S2.3, \S3.4] for
the details; there a proof is developed, and it
is precisely reproduced in [20]. By the way, it is possible to improve
$(4.2)$ into the form same as $(4.5)$ but with $M=N=\sqrt{D}$; 
see [62], which gives an improvement upon $(4.3)$. {\it See the second Addendum given
above\/}.
\smallskip
\item{65)} This should be compared with M.N. Huxley [25].
Another striking application is Iwaniec
[30], which stands for the hitherto best approximation to Gauss'
conjecture on the existence of primes of the form $n^2+1$.
\smallskip
\item{66)} See [61, \S4.2] as well as [9].
\vskip 1.5cm
\noindent
{\bf References}
\medskip
\item{[1]} N.C. Ankeny and H. Onishi. The general sieve. 
Acta Arith., {\bf 10} (1964/65), 31--62.
\item{[2]} M.B. Barban. The large sieve method and its application 
to number theory. Uspehi Mat.\ Nauk, {\bf 21} (1966), 51--102. (Russian)
\item{[3]} H. Bohr and E. Landau. Beitr\"age zur Theorie der 
Riemannschen Zetafunktion.  Math.\ Ann., {\bf 74} (1913), 3--30.
\item{[4]} E. Bombieri. On the large sieve. 
Mathematika, {\bf 12} (1965), 201--225.
\item{[5]} E. Bombieri. The asymptotic sieve. 
Rend.\ Accad.\ Naz., {\bf XL} (5) 1/2 (1975/76), 243--269.
\item{[6]} E. Bombieri. Le Grand Crible dan la  
Th\'eorie Analytique des Nombres (Second \'Edition). 
Ast\'erisque {\bf 18}, Paris 1987.
\item{[7]} E. Bombieri, J.B. Friedlander and H. Iwa\-niec. 
Primes in arithmetic progressions to large moduli. Acta Math., {\bf 156}
(1986), 203--251.
\item{[8]} R.W. Bruggeman. Fourier coefficients of cusp forms. 
Invent. math., {\bf 45} (1978), 1--18.
\item{[9]} R.W. Bruggeman and Y. Motohashi. 
A new approach to the spectral theory of the fourth moment of the Riemann
zeta-function. Crelle's Journal, in print.
\item{[10]} V. Brun. \"Uber das Goldbachsche 
Gesetz und die Anzahl der Primzahlpaare. 
Arch.\ Mat.\ Natur.\ B, {\bf 34},
no.\ 8, 1915. 
\item{[11]} V. Brun. La s\'erie ${1\over5}+{1\over7}
+{1\over11}+{1\over13}+{1\over17}
+{1\over19}+{1\over29}+{1\over31}+{1\over41}
+{1\over43}+{1\over59}+{1\over61}+\cdots$ 
o\'u les d\'enominateurs sont ``nombres premiers jumeaux'' 
est convergente
ou finie. Bull.\ Sci.\ Math., (2) {\bf 43} (1919), 124--128.
\item{[12]} V. Brun. Le crible d'Eratosth\`ene et le th\'eor\`eme de
Goldbach. Videnskaps.\ Skr.\ Mat.\ Natur.\ Kl.\ Kristiana, no.\ 3, 1920. 
\item{[13]} A.A. Buchstab. New improvements in 
the me\-thod of the sieve of Eratosthenes. Mat. Sbornik (N.S.), {\bf 4}(46)
(1938), 375-387. (Russian)
\item{[14]} A.A. Buchstab. A combinatorial 
intensification of the sieve of Eratosthenes. Uspehi Mat.\ Nauk, {\bf 22}
(1967), 199--226. (Russian)
\item{[15]} J.-R. Chen.  On the representation 
of a large even integer as the sum of a prime and the product of at most
two primes. Sci.\ Sinica, {\bf 16} (1973), 157--176. 
\item{[16]} J.-R. Chen. On the distribution 
of almost primes in an interval. Sci.\ Sinica, {\bf 18} (1975), 611--627.
\item{[17]} J.B. Friedlander and H. Iwaniec.  
On Bom\-bieri's asymptotic sieve. Ann.\ Scuola
 Norm.\ Sup.\ Pisa Cl.\ Sci., (4) {\bf 5} (1978), 719--756.
\item{[18]} P.X. Gallagher. A large sieve density estimate near $\sigma=1$.
Invent.\ math., {\bf 11} (1970), 329--339.
\item{[19]} S. Graham. Applications of sieve methods. Ph.D. Dissertation, Univ.\
of Michigan, 1977.
\item{[20]} G. Greaves. Sieves in Number Theory. 
Springer-Verlag, Berlin etc.\ 2001.
\item{[21]} G. Hal\'asz. \"Uber die Mittelwerte 
multiplikativer zahlentheoretischer Funktionen. Acta Math.\ Acad.\ Sci.\
Hungar., {\bf 19} (1968), 365--403.
\item{[22]} H. Halberstam and H.-E. Richert. Sieve Methods.  
Academic Press, London etc.\ 1974.
\item{[23]} G.H. Hardy and J.E. Littlewood. Some problems 
of ``partitio numerorum''; III:
 On the expression of a number as a sum of primes. Acta Math., {\bf 44}
(1922), 1--70.
\item{[24]} G. Hoheisel. Primzahlprobleme in der  Analysis. Sitz.\
Preuss.\ Akad.\ Wiss., {\bf 33} (1930), 3-11.
\item{[25]} M.N. Huxley. On the difference between consecutive primes.
Invent.\ math., {\bf 15}
 (1972), 164--170.
\item{[26]} A. Ivi\'c. The Riemann Zeta-Function.  Theory and
Applications. Dover Publ., Inc., Mineola, New York 2003.
\item{[27]} H. Iwaniec. On the error term in the linear sieve. 
Acta Arith., {\bf 19} (1971), 1--30.
\item{[28]} H. Iwaniec. A half dimensional sieve.  Acta Arith., {\bf
29} (1976), 69--95.
\item{[29]} H. Iwaniec. The sieve of Eratosthenes--Legendre. 
Ann.\ Scuola Norm.\ Sup.\ Pisa, (4) {\bf 4} (1977), 257--268.
\item{[30]} H. Iwaniec. Almost-primes represented 
by quadratic polynomials. Invent.\ math., {\bf 47} (1978), 171--188.
\item{[31]} H. Iwaniec. Sieve methods. Intern. 
Cong\-ress of Math.\ Proc., Helsinki 1978,  Acad.\ Sci.\ Fennica, Helsinki
1980, pp.\ 357--364.
\item{[32]} H. Iwaniec. Fourier coefficients of cusp 
forms and the Riemann zeta-function.
 Exp.\ No.\ {\bf 18}, S\'em.\ Th.\ Nombres, Univ.\ Bordeaux 1979/80.
\item{[33]} H. Iwaniec. Rosser's sieve. Acta Arith., 
{\bf 36} (1980), 171--202.
\item{[34]} H. Iwaniec. A new form of the error 
term in the linear sieve.
Acta Arith., {\bf 37} (1980), 307--320.
\item{[35]} H. Iwaniec. Rosser's sieve -- bilinear 
forms of the remainder terms -- some applications. 
In: Recent Progress in
Analytic Number Theory, Vol.1, (eds.\ H. Halberstam and C. Hooley),
Academic Press,
 London 1981, pp.\ 203--230.
\item{[36]} H. Iwaniec and M. Jutila. Primes in short 
intervals. Ark.\ Mat., {\bf 17} (1979), 167--176.
\item{[37]} W.B. Jurkat and H.-E. Richert. 
An improvement of Selberg's sieve method.\ I. Acta Arith., 
{\bf 11} (1965), 217--240.
\item{[38]} M. Jutila. On Linnik's constant. Math.\ Scand., {\bf 41} (1977),
45--62.
\item{[39]} P. Kuhn. Zur Viggo Brun'schen 
Siebme\-thode.\ I. Norske Vid.\
Selsk.\ Forh., Trondhjem, {\bf 14} (1941), 145--148.
\item{[40]} N.V. Kuznetsov. Petersson hypothesis 
for forms of weight zero and Linnik hypothesis. Preprint, Khabarovsk
Complex  Res.\ Inst., East Siberian Branch Acad.\ Sci. USSR, Khabarovsk
1977. (Russian)
\item{[41]} Yu.\ V. Linnik.  The large sieve. C.R.  Acad.\ Sci.\ URSS
(N.S.), {\bf 30} (1941), 292--294.
\item{[42]} Yu.\ V. Linnik. On the least prime 
in an arithmetic progression.\ I. The basic theorem. Rec.\ Math.\
(Sbornik), {\bf 15} (1944), 139--178; II. The Deuring--Heilbronn
phenomenon. ibid., 347--368.
\item{[43]} Yu.\ V. Linnik. Dispersion Method 
in Binary Additive Problems. Leningrad 1961. (Russian)
\item{[44]} Yu.\ V. Linnik. Additive problems and  eigenvalues of the
modular operators.  Proc.\ Intern.\ Cong.\ Math.\ Stockholm, 
1962, pp.\ 270--284.
\item{[45]} J. Merlin. Sur quelques th\'eor\`emes
d'Arithm\'etique et un \'enonc\'e qui les contient. 
C.R. Acad.\ Sci.\ Paris, {\bf 153} (1911), 516--518.
\item{[46]} H.L. Montgomery. Topics in Multiplicative Number 
Theory. Lect.\ Notes
in Math., {\bf 227}, Springer-Verlag, Berlin etc 1971.
\item{[47]} H.L. Montgomery. The analytic principle of 
the large sieve. Bull.\ Amer.\ Math.\ Soc., {\bf 84} (1978), 547--567.
\item{[48]} H.L. Montgomery and R.C. Vaughan.  The large sieve.
Mathematika, {\bf 20} (1973), 119--134.
\item{[49]} Y. Motohashi. On some improvements of the 
Brun--Titchmarsh theorem. J. Math.\ Soc.\ Japan, 
{\bf 26} (1974), 306--323.
\item {[50]} Y. Motohashi. On a density theorem of Linnik. Proc.\ Japan Acad., 
{\bf51} (1975), 815--817.
\item{[51]} Y. Motohashi. An induction principle 
for the generalization of Bombieri's
prime number theorem. Proc.\ Japan Acad., {\bf 52} (1976), 273--275.
\item {[52]} Y. Motohashi. On the Deuring--Heilbronn phenomenon.\ I. 
Proc.\ Japan Acad., {\bf 53} (1977), 1--2; II, ibid, 25--27.
\item {[53]} Y. Motohashi. On Gallagher's prime 
number theorem.  Proc.\ Japan  Acad., {\bf 53} (1977), 50--52.
\item {[54]} Y. Motohashi. A note on the large sieve.  Proc.\ Japan Acad., {\bf 53}
(1977), 17--19; II. ibid, 122--124; III. ibid, {\bf 55} (1979), 92--94;
IV. ibid, {\bf 56} (1980), 288--290.
\item {[55]} Y. Motohashi. On Vinogradov's zero-free region for the Riemann
zeta-function. Proc.\ Japan Acad., {\bf 54} (1978), 300--302.
\item {[56]} Y. Motohashi. On the zero-free region of Dirichlet's 
$L$-functions. Proc.\ Japan Acad., {\bf 54} (1978), 332--334.
\item{[57]} Y. Motohashi. Primes in arithmetic 
progressions. Invent.\ math., {\bf 44} (1978), 163--178.
\item{[58]} Y. Motohashi. On the linear sieve.  Proc.\ Japan Acad., 
{\bf 56} (1980), 285--287; II. ibid, 386--388.
\item {[59]} Y. Motohashi. An elementary proof of Vinogradov's 
zero-free region for the Riemann zeta-function. 
In: Recent Progress in Analytic Number
Theory, Vol.1, (eds.\ H. Halberstam et al.\ ), Academic Press, London 1981, pp.\
257--267.
\item{[60]} Y. Motohashi. Sieve Methods and Prime Number Theory.
Tata Inst.\ Fund.\ Res.\ Lect.\ Math.\ Phy., {\bf 72}, Bombay 1983.
\item{[61]} Y. Motohashi. Spectral Theory of the 
Riemann Zeta-Function. Cambridge
Univ.\ Press, Cambridge 1997.
\item{[62]} Y. Motohashi. On the error term in the 
Selberg sieve. In: Number Theory in Pro\-gress, {\bf 2}, Proc. Zakopane
Conf.\ 1997 (eds.\ K. Gy{\H o}ry et al.), W. de 
Gruyter, Berlin etc.\ 1999,
pp.\ 1053--1064.
\item{[63]} K. Muroi. Extraction of square roots in Babylonian
mathematics. Historia Sci., {\bf 9} (1999), 127--133.
\item{[64]} K. Prachar. Primzahlverteilung. Springer 
Verlag, Berlin etc.,
1957.
\item{[65]} A. R\'enyi. On the representation of an even 
number as the sum of a prime and an almost prime. 
Izv.\ Akad.\ Nauk SSSR
Ser.\ Mat., {\bf 12} (1948), 57--78. (Russian)
\item{[66]} A. Selberg. On the zeros of Riemann's 
zeta-function. Skr.\ Norske Vid.\
Akad. Oslo (1942), No.\ 1.
\item{[67]} A. Selberg. Contribution to the theory of 
the Riemann zeta-function. Arch.\ f\"or Math.\ og Naturv.\ B, {\bf 48}
(1946), No.\ 5.
\item{[68]} A. Selberg. On an elementary method in the 
theory of primes. Norske Vid.\ Selsk.\ Forh., Trondhjem, {\bf 19} (1947),
64--67. 
\item{[69]} A. Selberg. On elementary methods in prime 
number theory and their limitation\-s. 11 th.\ Skand.\ Math.\ Kongr.,
Trondhjem, (1949), 13--22. 
\item{[70]} A. Selberg. The general sieve-method and 
its place in prime number theory. Proc.\ Intern.\ Cong.\ Math.,  Cambridge,
Mass., {\bf 1} (1950), 286--292. 
\item{[71]} A. Selberg. On the estimation of Fourier 
coefficients of modular forms.
Proc. Symp.\ Pure Math., {\bf 8} (1965), 1--15.
\item{[72]} A. Selberg. Sieve methods. Proc.\ Symp.\ Pure 
Math., {\bf 20} (1971), 311--351. 
\item{[73]} A. Selberg. Remarks on sieves. Proc.\ 1972 Number 
Theory Conf., Boulder 1972, pp.\ 205--216. 
\item{[74]} A. Selberg. Remarks on multiplicative functions.
Springer Lect.\ Notes in Math., {\bf 626} (1977), pp.\ 232--241. 
\item{[75]} P. Tur\'an. Eine neue Methode in der Analysis 
und deren Anwendungen.
Akad. Kiad\'o, Budapest 1953.
\item{[76]} R.C. Vaughan. An elementary method in prime 
number theory. Acta Arith., {\bf 37} (1980), 111--115.
\item{[77]} A.I. Vinogradov. The density hypothesis for 
Dirichlet $L$-series. Izv.\ Akad.\ Nauk SSSR Ser.\ Mat., {\bf 29} (1965),
903--934; Corrigendum. ibid. {\bf 30} (1966), 719--720. (Russian)
\vskip 1cm
\noindent
{\font\small=cmr7\small
\baselineskip=8pt
Yoichi Motohashi
\smallskip\noindent
Department of Mathematics,
\par\noindent
Nihon University,
\par\noindent
Surugadai, Tokyo 101-8308, Japan
\par\noindent
ymoto@math.cst.nihon-u.ac.jp
\par\noindent
www.math.cst.nihon-u.ac.jp/$\scriptscriptstyle\sim$ymoto/
\vfill}

\bye